%% file: main.tex
\documentclass[11pt]{article}
\usepackage{amsmath}
\numberwithin{equation}{section}
\usepackage{amsfonts,amsthm,amssymb,bm,bbm}
\usepackage{graphicx}
\usepackage{hyperref,setspace,natbib}
\usepackage{xurl}
\usepackage{float}
\usepackage{placeins}
\usepackage{booktabs,multicol,multirow}
\usepackage{cleveref}
\usepackage{threeparttable}
\usepackage[svgnames]{xcolor}
\usepackage{mleftright} %
\usepackage[ruled,vlined]{algorithm2e}

\usepackage{enumitem}

\usepackage{subcaption} %
\usepackage{caption}
\usepackage{makecell}
\usepackage{tikz}
\usetikzlibrary{arrows.meta, positioning}

\usepackage{threeparttable}
\usepackage{tabularx}
\newcolumntype{Y}{>{\centering\arraybackslash}X}

\input{math_command}
\newtheorem{rmk}{Remark}
\begin{document}

\title{~~~ \bf \LARGE \\[-.75in] 
\parbox[c]{\textwidth}{\centering 
A Dynamic Inventory Rationing Model for Last-Mile Delivery of Medication in Liberia's Community Health Program
}%
}
\author{
Baris Ata%
\thanks{Booth School of Business, University of Chicago. \textbf{Email:} \url{baris.ata@chicagobooth.edu}}, 
Yufan Chen%
\thanks{Booth School of Business, University of Chicago. \textbf{Email:} \url{yufan.chen@chicagobooth.edu}}, and
Robert Montgomery%
\thanks{Rady School of Management, University of California San Diego. \textbf{Email:} \url{romontgomery@ucsd.edu}}, 
}
\date{
October 22, 2026
}
\maketitle

\thispagestyle{empty}
\centerline{\bf Abstract} \baselineskip 14pt \vskip 10pt
\noindent
\noindent
\noindent

We study the dynamic allocation of scarce medications in Liberia’s community health program. A central warehouse supplies a health facility and the Community Health Assistants (CHAs) it supervises. When a CHA runs out of medication, patients are referred to the facility, incurring travel costs if they complete the referral and remain untreated otherwise. Allocating inventory to CHAs allows more patients to receive treatment without traveling, while retaining inventory centrally preserves flexibility to respond to uncertain demand. The allocation problem must balance these considerations while accounting for fairness across locations.

We formulate this problem as a finite-horizon stochastic control problem whose objective combines the costs of referrals and untreated patients with a demand-weighted fairness utility. Because standard dynamic programming methods are computationally intractable in this setting, we develop a continuous-time approximation that involves controlling the drift rate of a reflected Brownian motion. We then adapt the deep splitting method, using deep neural networks to approximate the value function and its gradient and compute a dynamic allocation policy.

We calibrate the model using data on malaria cases among children in the Fangonda region of Lofa County, Liberia, and evaluate the proposed policy in the original discrete-time model. In the baseline setting, our policy reduces operating costs by about 6\% relative to the strongest benchmark, the fairness-only policy, while retaining approximately 98\% of its fairness utility. Additional experiments show that higher referral costs shift inventory toward CHAs and reduce local unmet demand, while greater demand variability increases our policy's performance over the fairness-only benchmark. These findings illustrate the value of dynamic policies. More broadly, our model and computational approach provide a framework for studying similar last-mile allocation problems in which scarce resources must be distributed dynamically across multiple service locations.


\clearpage

\setcounter{page}{1}
\pagenumbering{arabic}
\setlength{\baselineskip}{1.5\baselineskip}
\onehalfspacing
\setstretch{1.45}
\setcounter{equation}{0}
\setcounter{table}{0}
\setcounter{figure}{0}

\section{Introduction}
\label{sec:introduction}

Liberia's health system relies on Community Health Assistants (CHAs), who provide integrated community case management (iCCM) diagnosis and treatment during regular household visits. The widespread iCCM approach focuses on malaria, pneumonia, and diarrhea cases in children under the age of five. The CHA program is intended to extend the reach of established health facilities to serve populations that would otherwise face difficulties in accessing healthcare, and is most effective when CHAs are not only able to diagnose but also treat cases they encounter in the field. When CHAs diagnose cases but are unable to treat patients during a household visit, they must refer these patients to the nearest regional health facility. There are two principal causes for CHAs being unable to treat patients: either the patients require care that is beyond the scope of CHA training (usually when the patient is too sick or the case involves complicating factors), or the CHA does not have the necessary medications to treat a case he or she would have otherwise been able to treat in the field. The latter of these two causes results from a point-of-care stockout. While the set of referrals resulting from required care beyond the scope of CHA practice does not have a straightforward solution, the number of referrals resulting from point-of-care stockouts can be reduced through better inventory management strategies. 

This aim is a point of focus for our implementing partner, The International Rescue Committee's (IRC) Liberia country office, which identified the reduction of CHA health facility referrals as an important area for improvement within the community health system. These referrals can place a significant burden on families, often in the form of travel costs and lost wages. In most cases, a referral requires the child and caregiver to either walk or pay for transportation to the nearest health facility. In Lofa County specifically, where we source data for this project, the average distance between the communities served by CHAs and health facilities is greater than 5km. 

The health facilities in Liberia involved in the CHA program each oversee a set of CHAs serving communities within its region. Community Health Services Supervisors (CHSS) attached to these health facilities are responsible for the ongoing training and supervision of the CHAs. In this role, CHSSs are the main source of communication between CHAs and the health facilities, as well as a means for transporting medications to CHAs. The health facility managers must decide how much medication to retain at the health facility and how much to distribute among CHAs. Allocating more medication to CHAs should reduce the number of referrals, but increases the risk of stockouts at the health facility. Additionally, the distribution of more medication to all CHAs without accounting for demand patterns may leave some CHAs with excess, unused supplies. We model this setting to acount for the trade-offs between holding more inventory centrally as opposed to dispersing more medication among CHAs, while also accounting for demand pattern differences between CHA locations.

Focusing on malaria cases among children under the age of five, we formulate an inventory rationing problem within a two-echelon supply network and develop a dynamic allocation policy. The system consists of a central health facility and multiple geographically dispersed CHAs. We model demand at each location as a binomial random variable, which we approximate using a Gaussian distribution. To balance illness and travel costs with geographical equity considerations across the populations served by the CHAs and the health facility, we adopt an objective function that combines these costs with an $\alpha$-fair utility function. 

As a solution for this objective, we pursue a dynamic rationing policy for this problem. Standard dynamic programming methods are computationally intractable in this setting due to the curse of dimensionality. We therefore develop a continuous-time approximation that involves controlling the drift rate of a reflected Brownian motion (RBM) in the orthant. To solve this approximate model, we develop a computational method based on deep neural networks. Specifically, we build on the deep splitting method of \cite{beck2021deep} to approximate the value function and its gradient, enabling us to solve the problem in high dimensions.

We calibrate our model using data from the Fangonda region of Lofa County, Liberia, and evaluate the proposed policy in the original discrete-time model. Our policy outperforms the benchmark policies considered in this paper. Relative to the fairness-only policy, the strongest benchmark, our policy reduces operating costs by about 6\% while retaining approximately 98\% of the benchmark's fairness utility. These results show that dynamically adjusting allocations based on the system state can substantially improve the balance between operating costs and allocation fairness.

Our additional experiments examine how the proposed policy responds to changes in the operating environment. We find that increasing the weight placed on fairness yields higher fairness utility at the expense of somewhat higher operating costs. Higher referral costs shift inventory toward CHAs, reducing local unmet demand and supply--demand mismatches. Greater demand variability makes demand more difficult to meet but also increases our policy's performance over the fairness-only benchmark. Higher spatial demand correlation reduces mismatches across locations, while changes in HF demand and initial inventory prompt the policy to adjust the balance of shipments between the HF and CHAs. Together, these experiments show that the proposed policy responds systematically to changes in costs, demand, and resource availability.

Overall, our results suggest that dynamic allocation can improve the use of scarce medication inventory by accounting for current inventory levels, uncertain demand, referral costs, and fairness considerations rather than relying solely on fixed allocation shares. Although grounded in the Liberian context, our model and computational approach also provide a framework for studying similar last-mile allocation problems in which limited resources must be distributed dynamically across multiple service locations.

The remainder of this paper is organized as follows. Section~\ref{sec:literature_review} provides a literature review considering inventory rationing models, humanitarian and global health supply chains, and techniques for solving high-dimensional Brownian control models. Section~\ref{sec:model} presents the model using a single warehouse multiple retailer framework. Section~\ref{sec:finite_horizon} establishes the deep slitting algorithm used to solve the HJB Equation. Section~\ref{sec:computational_method} develops the deep splitting method we use to approximate the value function and its gradient. Section~\ref{sec:real_world_application} puts forth a series of practical benchmark policies by which to compare the performance of the optimal policy. Section~\ref{sec:computational_results} uses data from CHA cases, and treaments in Lofa county along with health facility inventory stock levels to test the performance of proposed policies. Section~\ref{sec:conclusion} concludes the paper with a summary of the findings.

\section{Literature Review}
\label{sec:literature_review}

At its core, the modeling approach in this paper addresses a fundamental question: how should a health facility and its community health assistants (CHAs) share a limited medication supply? Inherently, this is a problem of inventory rationing. The tension in our model is captured, in part, through competing costs: keeping stock centrally located increases flexibility, whereas dispersing it among CHAs brings treatment closer to patients. This dynamic is common in distribution systems, and we draw on existing literature to understand the tradeoffs between these inventory placement strategies. Given the global health context of our setting, we also incorporate fair allocation measures to consider when constructing and evaluating policies. Finally, our modeling structure and proposed solutions rely on stochastic control techniques. Broadly, we connect this allocation problem to literature on inventory rationing, distribution systems, fairness, global health logistics, and stochastic control.

Inventory rationing models help determine when to reserves scarce stock for demand classes with different shortage costs or service requirements. \citet{veinott1965optimal} studies nonstationary inventory with several demand classes, and \citet{topkis1968ordering} derives class- and time-dependent thresholds, under his model's assumptions, for deciding when to serve demand or preserve inventory. \citet{nahmias1981operating} analyze rationing under periodic and continuous review. These foundational models show how the value of reserving stock depends on demand priorities and replenishment opportunities; \citet{romeroSilva2026review} review the broader literature.

Later studies clarify how shortage assumptions shape allocation policies. \citet{deshpande2003threshold} coordinate replenishment and threshold rationing for two demand classes with backorders. \citet{arslan2007single} evaluate multiclass rationing through an equivalence with serial inventory systems under a specified threshold policy and backorder-clearing rule. \citet{ha1997rationing} establishes structural results for production and rationing with multiple demand classes and lost sales. In a health application, \citet{samii2012vaccines} study vaccine inventory rationing.

Our model also starts with a fixed supply, but links shortages across treatment locations through patient referrals. Some patients encountering a CHA stockout reach the facility; others remain untreated, and referred patients still depend on the facility having medication. Rationing therefore affects both immediate access and demand at the facility.

Our model can be seen as a one-warehouse, multiple-retailer (OWMR) model. OWMR models study how shared supply should be distributed across locations. \citet{eppen1979centralization} quantifies pooling benefits under normal demands and common linear costs, highlighting the role of demand correlation. \citet{eppenSchrage1981centralized} analyze centralized ordering through a stockless depot under an allocation-balance assumption. \citet{federgruen1984allocation,federgruen1984approximations} develop allocation policies and systematic approximations for dynamic multilocation systems.

\citet{jackson1988allocation} is particularly relevant because a warehouse allocates a shipment over a cycle before its next outside replenishment. Retaining stock allows later allocations to respond to additional demand information. \citet{jacksonMuckstadt1989risk} examine this risk-pooling benefit in a two-period model. Our setting shares this reserve-allocation problem, with treatment access and referral costs determining the consequences of inventory placement.

Central backup service also has precedents. \citet{anupindi2001framework} allow some unmet retailer demand to be served from central or other retailer inventories, while \citet{axsater2007heuristics} consider direct customer demand at the warehouse. Our model distinguishes deployable warehouse reserve from facility treatment stock, even when they are co-located. The reserve can supply any treatment location; facility stock serves local patients and referrals. This distinction connects the flexibility of central reserves to the value of local treatment under seasonal demand.

Our objective combines stockout-related costs with a demand-weighted alpha-fair utility of allocation rates relative to expected local demand rates. This utility rewards both the amount allocated and its distribution. Fairness objectives depend on what is being distributed and how recipients are compared. \citet{mo2002fair} develop the alpha-fair family for network resource allocation, and \citet{bertsimas2012efficiency} characterize its efficiency–fairness trade-offs. \citet{karsu2015inequity} review alternative formulations. In related applications, \citet{mccoy2014fairness} study fairness in allocating health-delivery transport capacity. 

Integrated community case management extends treatment of childhood malaria, pneumonia, and diarrhea to communities with limited facility access \citep{young2012joint}. A controlled before-and-after evaluation of a precursor program in Rivercess, Liberia, found increased treatment by qualified providers \citep{white2018liberia}. Reliable supplies are essential to this delivery model. \citet{chandani2012medicines} examine medicine availability among community health workers, and \citet{yadav2015supply} identifies broader weaknesses in health-product supply chains.

While the traditional cost of centralized inventory is increased delivery response time, modern literature also views decentralization as a risk mitigation tool. \cite{schmitt2015centralization} measured the benefits of decentralized inventory in reducing cost variance and lowering the probability of a stockout during a warehouse delay or failure. Additionally, with the growth of e-commerce, the OWMR framework has expanded to include warehouses as sources of direct order fulfillment, paralleling the direct-to-patient dynamics in our model. In our specific context, the response time cost of centralization translates directly to patient burden. Patients unable to receive treatment at decentralized CHAs incur the cost of traveling to a centralized facility. Consequently, some of them will not receive treatment. We rely on literature measuring the cost of failed referrals in community health settings, such as \cite{morris2023factors} and \cite{rasanathan2014here,rasanathan2014community}. In general, the limited studies on community-based referrals in sub-Saharan Africa indicate that pregnancy and newborn-related referrals have a much higher realized referral rate than iCCM conditions, with compliance for the latter falling below 10\% in some settings \citep{lal2018compliance}. In contrast, \citet{jarolimova2018completion} find substantially higher follow-through in another Ugandan program. \citet{nanyonjo2015referral} document referral-related expenses and caregivers' opportunity costs. These findings motivate including patient movement and access costs in inventory decisions, while recognizing that referral behavior varies across settings.

Several operations studies address closely related distribution problems. \citet{vledder2019zambia} evaluate a large supply-chain field experiment in Zambia. \citet{gallien2021zambia} use a field-informed, validated simulation to compare pharmaceutical inventory policies and examine availability and access equity. \citet{parvin2018malawi} provide a particularly close precedent: their Malawi model combines delayed shipments, clinic clustering, and dynamic transshipments under seasonal demand, with equity considerations. Our focus is on stockouts that move patients toward the supervising facility. Such referrals couple local availability with facility demand and patients' travel burdens.

We first formulate a discrete-time stochastic dynamic optimization model for inventory allocation, following the existing literature. However, this model is computationally intractable due to the curse of dimensionality. We then use a moment-matched Gaussian approximation of the demand process and adopt a continuous-time formulation. This yields a controlled reflected Brownian motion (RBM) model coupled with an exhaustible warehouse reserve. The reflection process tracks unmet demand and generates referrals to the facility when CHAs run out of stock. \citet{harrison1981reflected} establish foundational results for RBM, and \citet{harrison1988brownian} develops Brownian network models for dynamic control. Applications include make-to-stock scheduling \citep{wein1992dynamic}, joint capacity and inventory decisions \citep{bradley2002capacity}, and volunteer staffing in nonprofit gleaning operations \citep{ata2019gleaning}. Our construction in this paper is a demand-based approximation rather than a heavy-traffic limit theorem. \citet{ata2025drift} provide a closely related computational method for drift control of RBM under infinite-horizon objectives. 

We use deep neural-network based method to solve our dynamic allocation problem. Deep neural network methods approximate value functions and their marginal values from simulated paths, avoiding a full grid over inventory states. The deep backward stochastic differential equation (BSDE) methods of \citet{e2017deep} and \citet{han2018solving} train neural networks for the value function and its gradient. They propagate approximate values along paths and train against a known terminal condition, coupling parameters across time.

\citet{beck2021deep} develop deep splitting, which replaces joint training across all time steps with successive smaller regression problems. Our adaptation works backward from the terminal condition. At each step, the next-time value approximation is fixed while current values and allocation-relevant marginal values are fitted jointly. Training smaller problems eases computation.

We adopt the deep splitting to drift control of RBM. Our adaptation accounts for referral flows, stockout costs, and the warehouse resource constraint. It learns the marginal effect of sending another unit to a CHA instead of retaining it centrally. We apply this method to solve the inventory allocation problem for the Lofa County. That corresponds to a nine dimensional control problem. While this is remarkable by historical standards, we view it only as a proof-of-concept, and expect that our method can handle problems of much higher dimensions, e.g., in hundreds if not thousands of dimensions; see \cite{beck2021deep} and \cite{ata2025dynamic}. 

Our contribution combines dynamic medication rationing, stockout-induced referrals, unequal access costs, and an explicit preference for equitable allocation under limited supply. Adapting deep splitting to this problem allows us to examine how central reserve retention and local treatment access should be balanced over time.

\section{Model} 
\label{sec:model}

Our model assumes a single warehouse supplies a health facility (HF), indexed by $k=0$, and $K$ community health assistants (CHAs), indexed by $k=1, \dots, K$. The warehouse receives a shipment at the beginning of the planning horizon and no further shipments are received during that planning horizon $[0,T]$. CHAs provide treatment to patients within each of the $K$ service regions, where CHA $k$ serves region $k$ for $k=1, \dots, K$. The HF primarily serves patients in its vicinity, which we refer to as region $0$. The HF is also responsible for treating patients who travel to the HF after being referred from any of the $K$ CHAs.

A system manager is responsible for allocating the warehouse inventory to the HF and CHAs over the course of the planning horizon $[0,T]$. We divide this planning horizon into $N$ review periods of equal length $h=T/N$, and let $t_n=nh$ for $n=0,1,\dots,N$ denote the review epochs.  We use $n$ to index the review periods, where $n$ corresponds to the interval $[t_{n-1},t_n]$. In this way, the system manager is responsible for determining if and how much inventory should be distributed to the central HF and how much inventory should be distributed to the CHAs at the start of each period $n$. 

We use $\xi_k(n)$ to denote the demand for region $k$ in period $n$, where $k=0,1,\dots,K$ and $n=1,\dots,N$. Each demand $\xi_k(n)$ is modeled using a binomial distribution, where demand is assumed to be independent across locations and time, such that  
\begin{equation}
    \xi_0(n) \sim \mathrm{Binomial}(N_h,\theta(n)),\quad \text{and} \quad \xi_k(n) \sim \mathrm{Binomial}(N_c,\theta(n)), \ k=1,\dots,K.
    \label{eq:binomial_model}
\end{equation}
We define $N_h$ as the number of individuals assigned to the HF within its local catchment area, and $N_c$ as the number of individuals assigned to each CHA. We assume that each CHA serves the same number of individuals. Although in practice these population sizes can vary, the community health program is designed so that each CHA is assigned roughly the same number of families. We define $\theta(n)$ as the probability that a given individual in any location is diagnosed with malaria and determined to require a course of treatment during period $n$. We note that this probability is not dependent on whether an individual has previously had malaria, as previous malaria episodes do not confer reliable protection against reinfection or recurrent clinical malaria in children under five. This lack of immunity is observed in high-transmission settings where young children may experience frequent new infections and repeated symptomatic episodes despite substantial previous exposure \citep{tchaparian2016population,trape2024one}.

During each period, the CHA assigned to region $k$ (for $k=1,...,K$) supplies medication to the patients requiring treatment in region $k$ from the CHA's existing inventory. If demand during a period exceeds the existing inventory available for region $k$, the patients who are unable to receive treatment in region $k$ are referred to the HF. With probability $p_k$ (for $k=1,...,K$), a referred patient successfully follows the referral and travels to the HF. These patients who travel to the HF following a referral incur a referral cost of $r_k$ (for $k=1,...,K$). The referral cost captures the additional costs associated with travel and lost caregiver wages from needing to seek care at the health facility. Inversely, with probability $1-p_k$, a referred patient is unable to successfully follow the referral and does not receive treatment. In this case, the patient incurs a cost of $u$ (meaning the patient is untreated). The HF is thus responsible for treating both the patients from its own region (region $0$) who require treatment, as well as the patients who are referred by the $K$ CHAs and successfully travel to the HF. If the total demand at the HF exceeds the available inventory at the HF, the patients who are unable to receive treatment at the HF incur the cost $u$. This means that a patient who is referred from location $k$ to the HF, travels to the HF, and is then unable to receive treatment incurs a total cost of $r_k+u$.

In order to make an allocation decision for period $n+1$, the system manager reviews the system state at end of period $n$. We denote the system state at the end of the period $n$ by $(Z_{0}(n),Z_{1}(n),\dots,Z_{K}(n),Q(n))$, where $Z_{0}(n)$ denote the inventory level at the HF, $Z_{k}(n)$  denotes the inventory level at CHA $k$ for $k=1,2,\dots,K$, and $Q(n)$ denotes the warehouse inventory level. We assume that the HF and CHAs both begin with zero inventory at the start of the planning horizon, that is,
\begin{equation*}
    Z_{k}(0) = 0,\quad \text{for} \quad k=0,1,\dots,K.
\end{equation*}
At the beginning of each period, upon observing the system state, but before demand is realized, the system manager decides how much inventory to allocate from the warehouse to the HF and each CHA. Once inventories are allocated, they are received at each location instantaneously. Inventory levels are then updated, demand is realized,  patients are served, and costs are incurred. Then the period ends.

The system manager's allocation decision in period $n$ is defined by
\begin{equation}
    a(n) = (a_{0}(n),a_{1}(n),\dots,a_{K}(n))^{\top}, 
\end{equation}
where $a_{0}(n)$ denotes the amount of inventory allocated to the HF and $a_{k}(n)$ denotes the amount of inventory allocated to CHA $k$ for $k=1,2,\dots,K$. The system manager's control policy is then defined as $\{ a(n) \geq 0 : n \leq N \}$, and the cumulative allocation process is defined as
\begin{equation}
    A(0)=0, \quad \text{and} \quad  A(n) = \sum_{s=1}^{n} a(s),\  n \leq N.
\end{equation}
For feasibility, we require that 
\begin{equation}
    e \cdot a(n) \leq Q(n-1), \quad n \leq N,
    \label{eq:feasibility_condition_discrete}
\end{equation}
where $e$ is a $(K+1)$-dimensional vector of ones. This allocation constraint ensures that the amount of inventory allocated to the $K+1$ regions in period $n$ does not exceed the available warehouse inventory. Consequently, over the entire planning horizon, the total inventory allocated to the health facility and CHAs cannot exceed the starting balance of warehouse inventory, which is written as, 
\begin{equation}
    e \cdot A(N) \leq Q(0).
\end{equation}
We also require that
\begin{equation}
    a_k(n) \leq \bar{a}_k(n), \qquad n \leq N,
    \label{eq:shipment_bound_discrete}
\end{equation}
where $\bar{a}_k(n)$ is an upper bound on the allocation to location $k$ in period $n$. This constraint reflects the carrying capacity of a community health services supervisor (CHSS), who delivers medications to the CHAs. The warehouse inventory then evolves as
\begin{equation}
    Q(n) = Q(0) - e \cdot A(n), \quad n \leq N.
    \label{eq:warehouse_inventory}
\end{equation}

In order to define the inventory level $Z_{k}(n)$ for $k=1,2,\dots,K$ and $n \leq N$, we track the cumulative demand $D_k(n)$ at location $k$, for $k=0,1,\dots,K$, until the end of period $n$. To do this, we define
\begin{equation}
    D_k(n) = \sum_{s=1}^{n} \xi_k(s),
\end{equation}
and set $D_{k}(0) = 0$. We can also track the cumulative demand at CHA $k$, for $k=1,\dots,K$, at the end of period $n$ that is not met locally, or equivalently, the cumulative number of referrals. To do this, we let
\begin{equation}
L_{k}(n)=\max_{0\leq s \leq n}[D_k(s)-A_k(s)],\quad n=1,2,\dots,N,
\end{equation}
and set $L_{k}(0)=0$ for $k=1,\dots,K$.

\begin{figure}[!ht]
    \centering
    \includegraphics[width=0.9\textwidth]{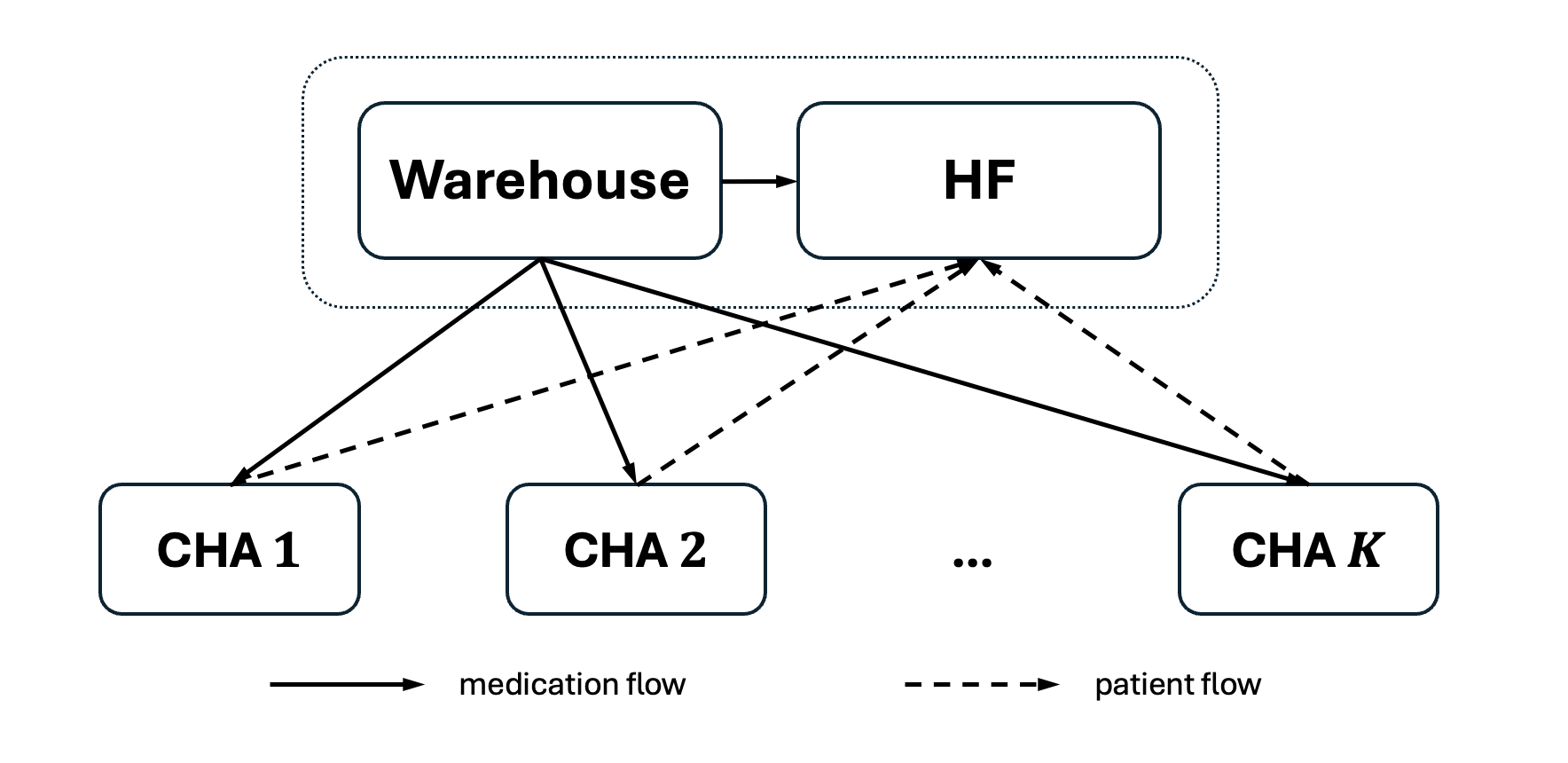}
    \caption{An illustration of the patient and medication flows within the model for a single Health Facility (HF) and $K$ Community Health Assistants (CHAs).}
    \label{fig:liberia_model}
\end{figure}

We can then define the inventory level $Z_{k}(n)$, for each CHA $k$ as, 
\begin{equation}
    Z_{k}(n) = Z_{k}(0) + A_{k}(n) - D_{k}(n) + L_{k}(n), \quad n \leq N.
    \label{eq:inventory_k}
\end{equation}
 Similarly, the inventory level of the HF evolves as
\begin{equation}
    Z_{0}(n) = Z_{0}(0) + A_{0}(n) - D_{0}(n) - \sum_{k=1}^{K} p_{k} L_{k}(n) + L_{0}(n), \quad n \leq N,
    \label{eq:inventory_0}
\end{equation}
where $L_{0}(n)$ denotes the cumulative unmet demand at the HF by the end of period $n$, and it is given as follows:
\begin{equation}
    L_{0}(n)=\max_{0\leq s \leq n} \left[D_0(s)+ \sum_{k=1}^{K} p_{k} L_{k}(s) - A_{0}(s)\right],
\end{equation}
with $L_0(0)=0$. The HF load combines its local demand with the expected number of complied referrals, $\sum_{k=1}^K p_kL_k(n)$. That is, in Equation~\eqref{eq:inventory_0}, we approximate the number of successfully completed referrals with its expectation for simplicity. 
Figure~\ref{fig:liberia_model} depicts the directional flows of medication and patients within the system. 

The system manager must then decide how much and where to allocate inventory in order to minimize costs. These costs include those previously discussed, namely: the cost associated with traveling to the health facility from location $k$ because of a referral ($r_k$), as well as the cost of a patient who goes untreated ($u$). Thus, the expected immediate cost associated with a stockout at location $k$ is written as
\begin{equation}
    c_{k} = \begin{cases}
        u, & k=0,\\
        (1-p_k)u + p_k r_k, & k =1,2,\dots,K.
    \end{cases}
    \label{eq:cost_definition}
\end{equation}

We also incorporate a fairness consideration for how medications are distributed across locations and periods relative to expected demand. In order to model this consideration, we first define
\begin{equation}
    \mu_{k}(n) = \mathbb{E}[\xi_k(n)]
\end{equation}
as the mean demand at location $k$ in period $n$. We then define the demand-normalized allocation shipment as
\begin{equation}
    f_k(n) = \frac{a_k(n)}{\mu_k(n)} \quad \text{for all } k,n,
\end{equation}
and let $f(n) = (f_{0}(n),f_{1}(n),\dots, f_{K}(n))^{\top}$ and $\mu(n) = (\mu_{0}(n),\mu_{1}(n),\dots, \mu_{K}(n))^{\top}$ for $n \leq N$. The fairness of an allocation decision is then measured using a demand-weighted $\alpha$-fair utility function $\Psi_{\alpha}$, defined as
\begin{equation}
    \Psi_{\alpha}(\rho;\mu) = \frac{1}{1-\alpha}
    \sum_{k=0}^{K} \mu_k \rho_k^{1-\alpha}\quad \text{for } \rho\in \R_{+}^{K+1},\; \mu\in \R_{++}^{K+1},
    \label{eq:fairness_utility}
\end{equation}
where $ \alpha \in (0,1)$, with larger values placing a stronger penalty on allocation imbalances.

Combining the fairness consideration with the stockout costs, the system manager's objective is to make the allocation decision dynamically over time to minimize
\begin{equation}
    \E\Big[c\cdot L(N) - \lambda\, \sum_{n=1}^{N} \Psi_{\alpha}(f(n);\mu(n)) \Big],
    \label{eq:objective_discrete}
\end{equation}
where $\lambda \geq 0$ is the weight assigned to the fairness utility, and the cumulative unmet demand at time $n$ is $L(n) = (L_{0}(n),L_{1}(n),\dots, L_{K}(n))^{\top}$ for $n \leq N$.

The formulation above is a Markov Decision Process (MDP). However, it is computationally prohibitive to solve because the standard MDP-based techniques suffer from the curse of dimensionality. With this in mind, we next describe a continuous-time approximation of the foregoing problem.

\paragraph{An approximating Brownian control problem (BCP).} The key simplification in our approximation is to view the system manager's problem in continuous time and allow for allocations decisions to likewise be made continuously over time. With a slight abuse of notation, we reuse the earlier notation and replace the cumulative allocation process with the following continuous-time process:
\begin{equation}
    A(t) = \int_{0}^{t} a(s) \dif s, \quad t\in[0,T].
\end{equation}
Here $a(t)$ denotes the allocation rate at time $t$, in contrast to the discrete allocation quantity $a(n)$ introduced above. We also approximate the cumulative demand process $D(\cdot)$ by a Gaussian process. Specifically, we assume
\begin{equation}
    D(t) = \int_{0}^{t} \mu(s) \dif s + \int_{0}^{t} \sigma(s) \dif B(s), \quad t\in[0,T],
\end{equation}
where $\mu(t) = (\mu_{0}(t),\mu_{1}(t),\dots, \mu_{K}(t))^{\top}$, $\Gamma(t) = \sigma(t)\sigma(t)^{\top}$ is a $(K+1)\times(K+1)$ positive-definite covariance matrix for $t \in [0,T]$, and $B(\cdot)$ is a $(K+1)$-dimensional standard Brownian motion. 
In addition, we set $\mu(t)$ and $\Gamma(t)$ as piecewise constant functions: for $n=1,2,\dots,N$ and $t\in(t_{n-1},t_n]$,
\[ \mu_0(t)= \frac{1}{h}N_h \,\theta(n), \quad \mu_k(t)= \frac{1}{h} N_c\,\theta(n), \; k=1,2,\dots,K,  \]
\[ \Gamma_{00}(t) = \frac{1}{h} N_h\,\theta(n)(1-\theta(n)), \quad \Gamma_{kk}(t)= \frac{1}{h} N_c\, \theta(n)(1-\theta(n)),\; k=1,\ldots,K, \]
where $h=t_n-t_{n-1}$ (for all $n=1,\ldots,N$) is the period length, and $\Gamma_{ij}(t)=0$ for $i\neq j$. This approximation matches the mean and variance of the Brownian demand increment over one review period with the discrete binomial demand.

The state process $(Z(t),Q(t))$ and the cumulative unmet demand process $L(t)$ are now viewed as continuous-time processes over $[0,T]$. Crucially, the system dynamics equations \eqref{eq:warehouse_inventory}, \eqref{eq:inventory_k} and \eqref{eq:inventory_0} remain valid with these modifications. Thus, we write for $t\in[0,T]$ that
\begin{equation}
    \begin{aligned}
        Q(t) &= Q(0) - e \cdot A(t),\\
        Z_{k}(t) &= Z_{k}(0) + A_{k}(t) - D_{k}(t) + L_{k}(t), \quad k=1,2,\dots,K,\\
        Z_{0}(t) &= Z_{0}(0) + A_{0}(t) - D_{0}(t) - \sum_{k=1}^{K} p_{k} L_{k}(t) + L_{0}(t).
    \end{aligned}
    \label{eq:state_process_continuous}
\end{equation}

The system manager's goal is to choose a non-anticipating allocation policy $\{a(t): t\in[0,T]\}$ to minimize total expected costs. For simplicity, we restrict attention to Markov controls, that is, we seek a measurable function $\varphi:[0,T]\times\R_{+}^{K+2} \to \R_{+}^{K+1}$ such that 
\begin{equation}
    a(t) = \varphi(t,Z(t),Q(t)), \quad t\in[0,T].
\end{equation}
Under such a control, the system state evolves as follows:
\begin{equation}
    Q(t) = Q(0) -  \int_0^t e \cdot \varphi(s,Z(s),Q(s)) \dif s, \quad t\in[0,T],
    \label{eq:continuous_state_Q}
\end{equation}
and
\begin{equation}
    Z(t) = Z(0) + \int_{0}^{t} \varphi(s,Z(s),Q(s)) \dif s - \int_{0}^{t} \mu(s) \dif s - \int_{0}^{t} \sigma(s) \dif B(s) + R L(t), \quad t\in[0,T],
    \label{eq:continuous_state_Z}
\end{equation}
where
\begin{equation}
R = \begin{pmatrix}
1 & - p_1 & -p_2 &\cdots & - p_{K}\\
0 & 1 & 0 & \cdots & 0 \\
0 & 0 & 1 & \cdots & 0 \\
\vdots & \vdots & \vdots & \ddots & \vdots \\
0 & 0 & 0 & \cdots & 1 
\end{pmatrix}.
\label{eq:reflection_matrix}
\end{equation}
As done earlier, we also require that the cumulative unmet demand process $L(t)$ is minimal in the following sense: For $k=0,1,\dots,K$,
\begin{equation}
    L_{k} \text{ is non-decreasing and increases only when } Z_{k}=0.
    \label{eq:Lk_def_word}
\end{equation}
This can be equivalently written as
\begin{equation}
    \int_{0}^{T} Z_{k}(t) \dif L_{k}(t) = 0, \quad k=0,1,\dots,K. 
    \label{eq:Lk_def_formula}
\end{equation}
Additionally, for feasibility we require that
\begin{equation}
    Q(t)\geq 0 \quad \text{and} \quad Z(t) \geq 0, \quad t\in[0,T].
\end{equation}

The system manager's objective can be written as
\begin{equation}
    \E\Big[c\cdot L(T) - \lambda \int_{0}^{T} \Psi_{\alpha}(\tilde{f}(t);\mu(t)) \dif t\Big],
\end{equation}
where $\tilde{f}(t) = (\tilde{f}_{0}(t),\tilde{f}_{1}(t),\dots, \tilde{f}_{K}(t))^{\top}$ with
\begin{equation}
    \tilde{f}_k(t) = \frac{\varphi_k(t,Z(t),Q(t))}{\mu_k(t)}, \quad k=0,1,\dots,K,\,t\in[0,T].
\end{equation}
We then denote the expected total cost during $[t,T]$ under control $\varphi$ starting in state $(z,q)$ at time $t$ as
\begin{equation}
    J(t,z,q;\varphi) = \E\Big[ c \cdot [L(T) - L(t)] - \lambda \int_{t}^{T} \Psi_{\alpha}(\tilde{f}(s);\mu(s)) \dif s \big| Z(t) = z,Q(t) = q \Big].
\end{equation} We also require that
\begin{equation}
    a_{k}(t) \leq \eta_k(t), \quad t\in[0,T],\ k=0,1,\dots,K,
\end{equation}
where $\eta_k(t)$ is an upper bound on the allocation rate to location $k$ at time $t$. It follows from \eqref{eq:shipment_bound_discrete} that
\[ \eta_k(t) = \frac{\bar{a}_k(n)}{h}, \quad t \in [t_{n-1}, t_n), \quad n \leq N. \]
Starting at time $t\in[0,T]$ in state $(z,q)$, the system manager seeks to choose a control $\varphi$ so as to minimize $J(t,z,q;\varphi)$.

Thus, we define the optimal value function $V$ as follows:
\begin{equation}
    V(t,z,q) = \inf J(t,z,q;\varphi), \quad t\in[0,T], \ (z,q)\in\R_{+}^{K+2},
\end{equation}
where the infimum is taken over feasible controls.

To facilitate the analysis to follow, we define the set of feasible allocation rates, denoted by $\mathcal{A}(t,q)$, as follows:
\begin{equation}
    \mathcal{A}(t,q) = \{a \geq 0: a_k \leq \eta_k(t) \ \text{for} \  k=0,1,\dots,K\} \text{ for } q>0,
\end{equation}
whereas when $q=0$ all allocation rates are set to zero, i.e., $\mathcal{A}(t,0)=\{0\}$.

\section{HJB Equation and the Proposed Policy} 
\label{sec:finite_horizon}
To characterize the optimal value function $V$ analytically, we present the associated Hamilton--Jacobi--Bellman (HJB) equation in this section. As a preliminary step, we let $f(t,z,q)$ be an arbitrary smooth function from $[0,T]\times \mathbb{R}_+^{K+1}\times \mathbb{R}_+$ to $\mathbb{R}$, and define the first-order differential operators $\nabla_z$, $\mathcal{R}$, and $\mathcal{D}$ as follows:
\begin{align*}
    \nabla_z f = \Big(\frac{\partial f}{\partial z_0}, \frac{\partial f}{\partial z_1},\dots,\frac{\partial f}{\partial z_K}\Big)^{\top}, \quad \mathcal{R} f = R^\top \nabla_z f,
\end{align*}
and
\begin{align*}
    \mathcal{D} f = \Big(\frac{\partial f}{\partial z_0} - \frac{\partial f}{\partial q}, \frac{\partial f}{\partial z_1} - \frac{\partial f}{\partial q}, \dots, \frac{\partial f}{\partial z_K} - \frac{\partial f}{\partial q}  \Big)^{\top}.
\end{align*}
Also, we define the second-order differential operator $\mathcal L_t$ for $t\in[0,T]$ via
\begin{align*}
    \mathcal{L}_{t} f = \frac{1}{2} \sum_{j=0}^{K} \sum_{k=0}^{K} \Gamma_{jk}(t) \frac{\partial^2 f}{\partial z_j \partial z_k} - \sum_{k=0}^{K} \mu_k(t) \frac{\partial f}{\partial z_k}.
\end{align*}
Then the HJB equation can be written as follows; see \cite{yong1999stochastic}, \cite{fleming2006controlled}, and \cite{ata2025drift}:
\begin{equation}
    \frac{\partial}{\partial t} V(t,z,q) + \mathcal{L}_t V(t,z,q) + \min_{a\in \mathcal{A}(t,q)}\Bigg\{a \cdot \mathcal{D} V(t,z,q) - \frac{\lambda}{1-\alpha} \sum_{k=0}^{K} \mu_k(t)\Big(\frac{a_k}{\mu_k(t)}\Big)^{1-\alpha} \Bigg\} = 0,
    \label{eq:hjb_finite}
\end{equation}
such that the following boundary and terminal conditions hold:
\begin{align}
& \mathcal{R}_k V(t,z,q) = - c_k,\quad  \text{if } z_k=0,\ k=0,1,\ldots,K,\quad q>0, \label{eq:hjb_boundary_condition_finite} \\
& V(T,z,q)=0, \quad z\ge 0,\quad q\ge 0. \label{eq:terminal_condition_finite}
\end{align}
The boundary conditions follow similarly to those in \cite{ata2025drift}, whereas the terminal condition follows from assigning zero cost to any remaining inventory.

Given the solution of the HJB equation, the optimal policy for the approximating Brownian control problem can be derived via the minimization term in Equation~\eqref{eq:hjb_finite}. This minimization problem is separable, and it leads to the following characterization of the optimal allocation rate $\varphi_k^*(t,z,q)$ for $k=0,1,\ldots,K$:
\begin{equation}
    \varphi_k^*(t,z,q)= 
    \begin{cases}
    0, & q=0,\\
    \eta_k(t), & q>0,\,\mathcal{D}_k V(t,z,q) \leq \lambda[\mu_k(t)/\eta_k(t)]^\alpha, \\
    \mu_k(t)\Big(\dfrac{\lambda}{\mathcal{D}_k V(t,z,q)}\Big)^{1/\alpha}, & q>0,\,\mathcal{D}_k V(t,z,q) > \lambda[\mu_k(t)/\eta_k(t)]^\alpha.
    \end{cases}
    \label{eq:optimal_policy_closed_form}
\end{equation}

These allocation rates are governed by two quantities: i) the remaining inventory $q$ at the warehouse; and ii) $\mathcal{D}_k V(t,z,q) = \frac{\partial V}{\partial z_k}(t,z,q) - \frac{\partial V}{\partial q}(t,z,q)$, which can be interpreted as the marginal opportunity cost of allocating one additional unit from the central inventory to location $k$: it measures the change in the continuation value when one unit is moved from the central stock $q$ to the local inventory $z_k$. The closed-form expression of the allocation rates in Equation~\eqref{eq:optimal_policy_closed_form} leads to the following insights:
\begin{itemize}
    \item As $\lambda\to 0$, the fairness term vanishes, and the policy becomes purely efficiency-driven. In this case, the optimal policy has a bang-bang form:
    \begin{align*}
    \varphi_k^*(t,z,q)=
    \begin{cases}
    \eta_k(t), & q>0 \text{ and } \mathcal D_k V(t,z,q)\le 0,\\
    0, & q>0 \text{ and } \mathcal D_k V(t,z,q)>0.
    \end{cases}
    \end{align*}
    \item When $\lambda>0$, fairness creates an additional incentive to allocate inventory, especially when the allocation level is low. As a result, the planner may allocate to location $k$ even when $\mathcal D_k V(t,z,q)>0$, provided that this marginal opportunity cost is not too large. Specifically, full allocation is optimal whenever 
    \begin{align*}
        \mathcal D_k V(t,z,q)\leq \lambda \Big[\frac{\mu_k(t)}{\eta_k(t)}\Big]^\alpha,
    \end{align*}
    while for larger values of $\mathcal{D}_k V(t,z,q)$, the optimal allocation rate decreases continuously according to
    \begin{align*}
    \varphi_k^*(t,z,q) = \mu_k(t)\Big(\dfrac{\lambda}{\mathcal{D}_k V(t,z,q)}\Big)^{1/\alpha}.
    \end{align*}
\end{itemize}

\paragraph{Proposed Policy.}
As will be described in Section~\ref{sec:computational_method}, our computational method uses a time grid $\{\tau_0,\tau_1,\ldots,\tau_M\}$ that is a refinement of $\{t_0,t_1,\ldots,t_N\}$. More specifically, we let
\[ \tau_m = \frac{m}{M}\,T, \quad m=0,1,\ldots,M, \]
where $r=M/N>1$ is an integer. Thus, we have
\[ t_n = \tau_{rn}. \]
Consequently, our computational method yields a neural network approximation $G_m$ of $\mathcal{D} V(\tau_m,\cdot,\cdot)$ for $m=0,1,\dots,M-1$. Substituting these into \eqref{eq:optimal_policy_closed_form} gives the proposed allocation rate functions $\varphi_k^*(\tau_m,\cdot,\cdot)$ at times $\tau_0,\tau_1,\ldots,\tau_M$.

To determine the proposed allocations to various CHAs and the HF in period $n$ at time $t_{n-1}$ $(n=1,2,\dots,N)$ given the system state $(Z(n-1),Q(n-1))$, we proceed as follows: We envision a virtual system that starts in state $(Z(n-1),Q(n-1))$ at time $t_{n-1}=\tau_{r(n-1)}$ and evolves through times $\tau_{r(n-1)},\tau_{r(n-1)+1},\ldots,\tau_{rn}$. Crucially, we assume there is no uncertainty, i.e., we consider a fluid model, and the functions $\varphi_k^*(\tau_m,\cdot,\cdot)$ $(m=r(n-1),r(n-1)+1,\ldots,rn-1)$ are used for making the allocation decisions at these times for the virtual system. That is, we let
\[ \hat{Z}(\tau_{r(n-1)})=Z(n-1), \quad \hat{Q}(\tau_{r(n-1)})=Q(n-1), \quad \hat{a}(\tau_{r(n-1)}) = \varphi_k^*( \tau_{r(n-1)},\hat{Z}(\tau_{r(n-1)}),\hat{Q}(\tau_{r(n-1)}) ), \]
and
\begin{align*}
    \hat{Z}(\tau_{m+1}) &= \hat{Z}(\tau_m) + [\hat{a}(\tau_m)-\mu(\tau_m)] (\tau_{m+1}-\tau_m) + R[\hat{L}(\tau_{m+1})-\hat{L}(\tau_m)], \\
     \hat{Q}(\tau_{m+1}) &= \hat{Q}(\tau_m) - e\cdot\varphi( \tau_m,\hat{Z}(\tau_m),\hat{Q}(\tau_m) ) (\tau_{m+1}-\tau_m),
\end{align*}
where $(\hat{Z},\hat{Q})$ and $\hat{a}$ denote the state process and the allocation decisions, respectively, in the virtual system. Additionally, $\hat{L}$ process is derived as described above and it corresponds to the minimal pushing process at the boundary that keeps $\hat{Z}$ nonnegative. Then, we set the tentative allocation for the HF and the CHAs in period $n$ by aggregating the allocations made in the virtual system as follows:
\[ \tilde{a}(n) = \sum_{m=r(n-1)}^{rn-1} \hat{a}(\tau_m)(\tau_{m+1}-\tau_m) = \frac{T}{M}\sum_{m=r(n-1)}^{rn-1} \hat{a}(\tau_m) , \quad n=1,\ldots,N. \]

If the total tentative allocation exceeds the available warehouse inventory, we allocate the remaining inventory in proportion to the tentative allocation vector. Thus, after enforcing the weekly shipment bounds, the actual allocation for $k=0,1,\dots,K$ is given by
\begin{equation}
a_k(n) =
\begin{cases}
\tilde{a}_k(n), & e \cdot \tilde{a}(n) \leq Q(n-1), \\
\tilde{a}_k(n) \cdot \dfrac{Q(n-1)}{e \cdot \tilde{a}(n)} , & e\cdot \tilde{a}(n) > Q(n-1).
\end{cases}
\label{eq:hjb_discrete_feasible_allocation}
\end{equation}

\section{Computational Method} 
\label{sec:computational_method}
Our computational method builds on the deep splitting method proposed by \cite{beck2021deep}. We first specify a \emph{reference policy}, denoted by $\tilde{\varphi}$, which is used to generate sample paths of the system state during neural-network training. In the reported training, reference paths come from fairness-only, static, and fluid-optimal policies (see Appendix~\ref{sec:implementation_details}). Intuitively, the reference policy should be chosen so that the resulting state process, which we refer to as the reference process, visits regions of the state space that the system is likely to occupy under an optimal policy. The corresponding reference processes $\tilde{Z}$ and $\tilde{Q}$ are then given as
\begin{equation}
\begin{aligned}
\tilde{Z}(t) &= \tilde{Z}(0) + \int_{0}^{t} [\tilde{\varphi}(s,\tilde{Z}(s),\tilde{Q}(s)) - \mu(s)] \dif s + W(t) + R \tilde{L}(t), \\
\tilde{Q}(t) &= \tilde{Q}(0) - \int_{0}^{t} e \cdot \tilde{\varphi}(s,\tilde{Z}(s),\tilde{Q}(s)) \dif s.
\end{aligned}
\label{eq:reference_rbm_finite}
\end{equation}
Recall that $\tilde{L}$ is minimal in the sense that it is nondecreasing and $\tilde{L}_{k}$ increases only when $\tilde{Z}_{k}=0$, cf. Equations~\eqref{eq:Lk_def_word}--\eqref{eq:Lk_def_formula}.  The sample paths of the reference processes serve as the training data for our computational method.

Next, we derive a stochastic identity based on the HJB equation~\eqref{eq:hjb_finite} by evaluating it at $(t, \tilde{Z}(t), \tilde{Q}(t))$ 
for $t\in[0, T]$ and by adding and subtracting the following term:
\begin{equation}
\begin{aligned}
\tilde{\varphi}(t,\tilde{Z}(t),\tilde{Q}(t)) \cdot \mathcal{D} V(t,\tilde{Z}(t),\tilde{Q}(t)),
\end{aligned}
\label{eq:stochastic_identity}
\end{equation}
which yields the stochastic identity:
\begin{equation}
    \begin{aligned}
        & \frac{\partial}{\partial t} V(t,\tilde{Z}(t),\tilde{Q}(t)) + \mathcal{L}_t V(t,\tilde{Z}(t),\tilde{Q}(t)) \\
        & \qquad \qquad + \tilde{\varphi}(t,\tilde{Z}(t),\tilde{Q}(t)) \cdot \mathcal{D} V(t,\tilde{Z}(t),\tilde{Q}(t)) + F(t,\tilde{Z}(t),\tilde{Q}(t),\mathcal{D}V(t,\tilde{Z}(t),\tilde{Q}(t))) = 0,
    \end{aligned}
    \label{eq:hjb_finite_stochastic}
\end{equation}
where for any $t\in[0,T]$, $z\in\R^{K+1}_{+}$, $q\in\R_{+}$ and $x\in \R^{K+1}$, the auxiliary function $F$ is defined as
\begin{equation}
\begin{aligned}
F(t,z,q,x) &= - \tilde{\varphi}(t,z,q) \cdot x + \min_{a\in \mathcal{A}(t,q)} \Bigg\{a \cdot x  - \frac{\lambda}{1-\alpha} \sum_{k=0}^{K} \mu_k(t) \Big(\frac{a_k}{\mu_k(t)}\Big)^{1-\alpha} \Bigg\}.
\end{aligned}
\label{eq:F_def}
\end{equation}

We proceed by fixing a partition $0 = \tau_0 < \tau_1< \dots < \tau_M = T$ of the finite horizon $[0,T]$ with $\Delta \tau_m = \tau_{m+1} - \tau_{m}$ for $m=0,1,\dots,M-1$. Then, given a reference policy $\tilde{\varphi}$, we generate the discretized state-pair samples, $(\tilde{Z}(\tau_m), \tilde{Q}(\tau_m))$ and $(\tilde{Z}(\tau_{m+1}), \tilde{Q}(\tau_{m+1}), \Delta\tilde{L}(\tau_{m+1}))$ with
\begin{align*}
    \Delta\tilde{L}(\tau_{m+1}) = \tilde{L}(\tau_{m+1}) - \tilde{L}(\tau_m),\quad m=0,1,\dots,M-1.
\end{align*}

To this end, given a reference policy $\tilde{\varphi}$, we first simulate discretized paths of the state process in Equation~\eqref{eq:reference_rbm_finite} with randomized initial states, as described in Remark~\ref{rmk:ref_path_perturbation}. At each training time $\tau_m$, the states reached by these reference paths form a set of states values at time $\tau_m$, denoted by $\mathcal{S}_m$. To simulate the reference process on the time interval $[\tau_m,\tau_{m+1}]$, we first sample $(\tilde{Z}(\tau_m), \tilde{Q}(\tau_m))$ uniformly from the set $\mathcal{S}_m$, and then simulate $(\tilde{Z}(\tau_{m+1}), \tilde{Q}(\tau_{m+1}), \Delta\tilde{L}(\tau_{m+1}))$ by applying a Euler--Maruyama discretization of Equation~\eqref{eq:reference_rbm_finite}. The sampling procedure is summarized in Algorithm~\ref{alg:sampling}, and the Euler--Maruyama step is given in Algorithm~\ref{alg:transition_kernel}.

\begin{rmk}[Reference-path sampling with initial perturbations] 
\label{rmk:ref_path_perturbation}
We generate the set $\mathcal{S}_m$ for $m=0,1,\dots,M-1$ from an ensemble of reference paths with randomized initial conditions. Let $Q_0$ denote the physical initial warehouse inventory. For 70\% of the reference paths, we use the physical initial condition
\[ \tilde{Z}(0)=0,\quad \tilde{Q}(0)=Q_0.\]
For the remaining 30\%, we perturb the initial condition by setting
\[ \tilde{Q}(0) = Q_0 + B_q, \quad \text{where} \quad B_q\sim \mathrm{Uniform}(-0.05Q_0,0.05Q_0),\]
and
\[ \tilde{Z}_k(0) \sim \mathrm{Uniform}(0,0.0125Q_0), \quad k=0,\dots,K,\]
independently across locations. The constants $70\%$, 0.05, 0.0125 should be viewed as tuning parameters and are chosen to increase exploration.
\end{rmk}

\begin{algorithm}[!ht]
\caption{Sampling$(m,B)$}
\label{alg:sampling}
\small
\KwIn{Time index $m$; batch size $B$; reference policy $\tilde{\varphi}$.}
\KwOut{ State-pair samples $\{(Z_m^{(b)},Q_m^{(b)}),(Z_{m+1}^{(b)},Q_{m+1}^{(b)},\Delta L_{m+1}^{(b)})\}_{b=1}^{B}$. }
\BlankLine
\For{$b \leftarrow 1$ \KwTo $B$}{
Sample one joint state vector $(Z_m^{(b)},Q_m^{(b)})$ uniformly from the set $\mathcal{S}_m$\;
Compute the reference allocation rate $ a_m^{(b)} = \tilde{\varphi}(\tau_m,Z_m^{(b)},Q_m^{(b)})$\;
Project $a_m^{(b)}$ onto the one-step feasible set so that $\sum_{k=0}^{K} a_{m,k}^{(b)}\Delta \tau_m \leq Q_m^{(b)}$\;
Generate $\Delta W_m^{(b)} \sim N(0,\Gamma(\tau_m)\Delta \tau_m)$\;
Set $ (Z_{m+1}^{(b)},Q_{m+1}^{(b)},\Delta L_{m+1}^{(b)}) = \textsc{StepTransition}(Z_m^{(b)},Q_m^{(b)},a_m^{(b)},\Delta W_m^{(b)},\Delta \tau_m)$\;
}
\Return{ 
    $\{(Z_m^{(b)},Q_m^{(b)}),(Z_{m+1}^{(b)},Q_{m+1}^{(b)},\Delta L_{m+1}^{(b)})\}_{b=1}^{B}$. 
    }
\end{algorithm}

\begin{algorithm}[!ht]
    \caption{StepTransition$(Z_m,Q_m,a_m,\Delta W_m,\Delta \tau_m)$}
    \label{alg:transition_kernel}
    \small
    \KwIn{Current state $(Z_m,Q_m)$; allocation rate $a_m$; Brownian increment $\Delta W_m$; step size $\Delta \tau_m$.}
    \KwOut{Next state $(Z_{m+1},Q_{m+1})$; regulator increment $\Delta L_{m+1}$.}
    \BlankLine
    Compute the unconstrained inventory vector $ X_{m+1} = Z_m + (a_m - \mu(\tau_m))\Delta \tau_m + \Delta W_m$\;
    For each CHA $k=1,\dots,K$, set $ \Delta L_{m+1,k} = \max\{0,-X_{m+1,k}\}$\;
    Set $ \Delta L_{m+1,0} = \max\left\{0, - \left(X_{m+1,0} - \sum_{k=1}^{K}p_k\Delta L_{m+1,k} \right) \right\}$\;
    Update the HF inventory: $Z_{m+1,0} = X_{m+1,0} + \Delta L_{m+1,0} - \sum_{k=1}^{K}p_k\Delta L_{m+1,k}$\;
    Update CHA inventories: $ Z_{m+1,k} = X_{m+1,k} + \Delta L_{m+1,k}$ for $k=1,\dots,K$\;
    Update warehouse inventory: $ Q_{m+1} = Q_m - \sum_{k=0}^{K}a_{m,k}\Delta \tau_m$\; \Return{$(Z_{m+1},Q_{m+1},\Delta L_{m+1})$.
    }
\end{algorithm}

In Equation~\eqref{eq:hjb_finite_stochastic}, the term involving the auxiliary function $F$ is nonlinear, whereas the remaining terms constitute a linear PDE if that term is ignored. The idea behind the deep splitting method is to approximate this nonlinear term in Equation~\eqref{eq:hjb_finite_stochastic} on each subinterval $[\tau_m, \tau_{m+1}]$ with a linear term, thereby obtaining a linear PDE that approximates the original PDE \eqref{eq:hjb_finite_stochastic} on that subinterval. One can then derive a Feynman-Kac representation for this linear PDE and use it to approximate the solution of the original PDE over the subinterval. This leads to the following, for $m=0,1,\dots,M-1$:
\begin{equation}
    \begin{aligned}
    & \E \Bigg[V(\tau_{m+1}, \tilde{Z}(\tau_{m+1}),\tilde{Q}(\tau_{m+1})) + c \cdot \Delta \tilde{L}(\tau_{m+1}) \Big| \tilde{Z}(\tau_m), \tilde{Q}(\tau_m) \Bigg] \\
    & \qquad = V(\tau_m, \tilde{Z}(\tau_m),\tilde{Q}(\tau_m)) - F(\tau_m, \tilde{Z}(\tau_m), \tilde{Q}(\tau_m), \mathcal{D} V(\tau_m, \tilde{Z}(\tau_m),\tilde{Q}(\tau_m))) \Delta \tau_m.
    \end{aligned}
    \label{eq:deep_splitting_relationship_approximation}
\end{equation}
See Appendix~\ref{appendix:derivation_of_deep_splitting_algorithm} for a formal derivation.

Crucially, it follows from Equation~\eqref{eq:deep_splitting_relationship_approximation} that one can focus on approximating $V$ and $\mathcal{D} V$ instead of $V$ and its gradient $\nabla V$; also see Equation~\eqref{eq:optimal_policy_closed_form}. Thus, for each $m=0,1,\dots,M-1$, we approximate the value function $V(\tau_m,z,q)$ with a neural network $H^{(m)}_{\omega_{m}}(z,q)$ and the function $\mathcal{D} V(\tau_m, z, q)$ with a neural network $G^{(m)}_{\nu_{m}}(z,q)$, where $\omega_{m}$ and $\nu_{m}$ are parameter vectors. The deep splitting method proceeds backward in time, starting with $m=M$. Using the terminal condition \eqref{eq:terminal_condition_finite}, we set $V(\tau_M,\tilde{Z}(\tau_{M}),\tilde{Q}(\tau_{M})) = 0$.

Because the conditional expectation in Equation~\eqref{eq:deep_splitting_relationship_approximation} can be interpreted as a projection, it can be characterized as the minimizer of a suitably chosen $L^2$ norm. Thus, to determine the optimal neural-network parameters, $\omega^{*}_{M-1}$ and $\nu^{*}_{M-1}$, we define the loss function as follows:
\begin{equation}
\begin{aligned}
& \mathcal{J}_{M-1}(\omega_{M-1},\nu_{M-1}) = \E\Bigg[\Big( - H^{(M-1)}_{\omega_{M-1}}(\tilde{Z}(\tau_{M-1}), \tilde{Q}(\tau_{M-1})) + c \cdot \Delta \tilde{L}(\tau_{M}) \\
& \qquad \qquad \qquad \qquad + F\Big( \tau_{M-1}, \tilde{Z}(\tau_{M-1}), \tilde{Q}(\tau_{M-1}), G^{(M-1)}_{\nu_{M-1}}(\tilde{Z}(\tau_{M-1}), \tilde{Q}(\tau_{M-1}))\Big) \Delta \tau_{M-1} \Big)^2 \Bigg]
\end{aligned}
\end{equation} 
where the expectation taken over $(\tilde Z(\tau_{M-1}),\tilde Q(\tau_{M-1}))$ and $(\tilde Z(\tau_{M}),\tilde Q(\tau_{M}),\Delta \tilde L(\tau_{M}))$, and we set $(\omega^{*}_{M-1}, \nu^{*}_{M-1}) \in \argmin \mathcal{J}_{M-1}(\omega_{M-1},\nu_{M-1})$.

Proceeding similarly for $m=M-2,\dots,1,0$, we first define the loss functions recursively as
\begin{equation}
\begin{aligned}
\mathcal{J}_{m}(\omega_{m},\nu_{m}) &= \E\Bigg[\Big( H^{(m+1)}_{\omega^*_{m+1}}(\tilde{Z}(\tau_{m+1}),\tilde{Q}(\tau_{m+1})) - H^{(m)}_{\omega_{m}}(\tilde{Z}(\tau_{m}),\tilde{Q}(\tau_{m})) + c \cdot \Delta \tilde{L}(\tau_{m+1}) \\
& \qquad \qquad \qquad \qquad \qquad + F\Big( \tau_{m}, \tilde{Z}(\tau_{m}), \tilde{Q}(\tau_{m}), G^{(m)}_{\nu_{m}}(\tilde{Z}(\tau_{m}),\tilde{Q}(\tau_{m})) \Big) \Delta \tau_{m} \Big)^2\Bigg]
\end{aligned}
\end{equation}
where the expectation taken over $(\tilde Z(\tau_{m}),\tilde Q(\tau_{m}))$ and $(\tilde Z(\tau_{m+1}),\tilde Q(\tau_{m+1}),\Delta \tilde L(\tau_{m+1}))$, and we set $(\omega^{*}_{m}, \nu^{*}_{m}) \in \argmin \mathcal{J}_{m}(\omega_{m},\nu_{m})$.

\begin{algorithm}[!ht]
    \caption{Deep Splitting Training via Backward Induction}
    \label{alg:training} \small
    \KwIn{
        Time grid $\{\tau_m\}_{m=0}^M$; batch size $B$; reference policy $\tilde{\varphi}$; optimization solver.\\ 
        Neural networks: value nets $H^{(m)}_{\omega_m}(z,q)$ and gradient-difference nets $G^{(m)}_{\nu_m}(z,q)$. }
    \KwOut{
        Trained value nets $H^{(m)}_{\omega_m^*}$, gradient-difference nets $G^{(m)}_{\nu_m^*}$, and parameters $\{\omega_m^*,\nu_m^*\}_{m=0}^{M-1}$.
    }
    \BlankLine
    Generate $B$ samples $\{(Z_{M-1}^{(b)},Q_{M-1}^{(b)}),(Z_M^{(b)},Q_M^{(b)},\Delta L_M^{(b)})\}_{b=1}^B$ by Algorithm~\ref{alg:sampling} over $[\tau_{M-1},\tau_M]$\;
    Initialize neural network weights $(\omega_{M-1},\nu_{M-1})$\;
    Compute the minimizer of the empirical loss $\hat{\mathcal{J}}_{M-1}(\omega_{M-1},\nu_{M-1})$ in Equation~\eqref{eq:empirical_loss_terminal}\;
    Set $(\omega_{M-1}^*,\nu_{M-1}^*)\in \argmin_{\omega_{M-1},\nu_{M-1}} \hat{\mathcal{J}}_{M-1}(\omega_{M-1},\nu_{M-1})$\;
    
    \For{$m=M-2,M-3,\ldots,0$}{
        Generate $B$ samples $\{(Z_m^{(b)},Q_m^{(b)}),(Z_{m+1}^{(b)},Q_{m+1}^{(b)},\Delta L_{m+1}^{(b)})\}_{b=1}^B$ by Algorithm~\ref{alg:sampling} over $[\tau_m,\tau_{m+1}]$\;
        Initialize neural network weights $(\omega_m,\nu_m)$\;
        Compute the minimizer of the empirical loss $\hat{\mathcal{J}}_m(\omega_m,\nu_m)$ in Equation~\eqref{eq:empirical_loss_recursive}\;
        Set $(\omega_m^*,\nu_m^*)\in \argmin_{\omega_m,\nu_m} \hat{\mathcal{J}}_m(\omega_m,\nu_m)$\;
    }
    
    \Return{$H^{(m)}_{\omega_m^*}$ and $G^{(m)}_{\nu_m^*}$ for $m=0,1,\dots,M-1$\;}
\end{algorithm}

We summarize our computational method in Algorithm~\ref{alg:training}, where we replace the population loss functions above with their empirical counterparts. More specifically, for each time interval $[\tau_m,\tau_{m+1}]$, we generate $B$ independent state-pair samples $\{(\tilde Z_m^{(b)},\tilde Q_m^{(b)}), (\tilde Z_{m+1}^{(b)},\tilde Q_{m+1}^{(b)},\Delta \tilde L_{m+1}^{(b)})\}_{b=1}^{B}$, where $(\tilde Z_m^{(b)},\tilde Q_m^{(b)})$ denotes the sampled state at time $\tau_m$, and $(\tilde Z_{m+1}^{(b)},\tilde Q_{m+1}^{(b)},\Delta \tilde L_{m+1}^{(b)})$ is generated from the Euler--Maruyama step under the reference policy $\tilde{\varphi}$; see Algorithm~\ref{alg:transition_kernel}. The expectation in $\mathcal{J}_m$ is then approximated by the sample average over these $B$ state pairs. For the terminal interval, this gives
\begin{equation}
\begin{aligned}
\hat{\mathcal{J}}_{M-1}(\omega_{M-1},\nu_{M-1}) &= \frac{1}{B}\sum_{b=1}^{B} \Big[ - H^{(M-1)}_{\omega_{M-1}}(\tilde Z_{M-1}^{(b)},\tilde Q_{M-1}^{(b)}) +  c\cdot \Delta \tilde L_{M}^{(b)}  \\
&\qquad + F\Big( \tau_{M-1}, \tilde Z_{M-1}^{(b)}, \tilde Q_{M-1}^{(b)},G^{(M-1)}_{\nu_{M-1}}(\tilde Z_{M-1}^{(b)},\tilde Q_{M-1}^{(b)})\Big) \Delta \tau_{M-1}\Big]^2 .
\end{aligned}
\label{eq:empirical_loss_terminal}  
\end{equation}
For $m=M-2,\dots,0$, the empirical loss is defined recursively as
\begin{equation}
\begin{aligned}
\hat{\mathcal{J}}_m(\omega_m,\nu_m) &= \frac{1}{B}\sum_{b=1}^{B}\Big[ H^{(m+1)}_{\omega_{m+1}^*}(\tilde Z_{m+1}^{(b)},\tilde Q_{m+1}^{(b)})- H^{(m)}_{\omega_m}(\tilde Z_m^{(b)},\tilde Q_m^{(b)}) + c\cdot \Delta \tilde L_{m+1}^{(b)}\\
&\qquad \qquad \qquad \qquad \qquad + F\Big(\tau_m,\tilde Z_m^{(b)},\tilde Q_m^{(b)},G^{(m)}_{\nu_m}(\tilde Z_m^{(b)},\tilde Q_m^{(b)})\Big)\Delta \tau_m \Big]^2 .
\end{aligned}
\label{eq:empirical_loss_recursive}
\end{equation}
The neural network parameters are obtained by minimizing these empirical losses backward in time.

As mentioned at the end of Section~\ref{sec:finite_horizon}, we choose the computation grid $\{\tau_0,\tau_1,\ldots,\tau_M\}$ as a refinement of $\{t_0,t_1,\ldots,t_N\}$. Thus, we learn the neural network approximation on that finer grid. Namely, we set 
\[ \tau_m = \frac{m}{M}\,T, \quad m=0,1,\ldots,M, \] 
where $r=M/N>1$ is an integer. Thus, we have $t_n=\tau_{rn}$.

\section{Data, Test Problems and Benchmark Policies}
\label{sec:real_world_application}

In this section, we calibrate our model using data from the Fangonda region of Lofa County, Liberia. As is typical of data collected in community health contexts, the dataset used for these purposes needs to be supplemented with secondary data sources. Using the information available, we attempt to capture the cost tradeoffs in the Fangonda and similar settings, while acknowledging its limitations. Furthermore, the model is intended to allow for insights that are broadly applicable to other contexts.

In what follows, we first describe the Fangonda data and explain how they inform the selected model parameters for our numerical experiments. We then define the test problem, introduce the benchmark policies used for comparison, and specify additional experiments.

\subsection{Data Description and Test Problem}
\label{sec:data_description}

The dataset includes records for seven CHAs associated with the HF. For each CHA, we observe weekly records of household visits, cases identified for pneumonia, malaria, and diarrhea, disease-specific treatments, and referrals issued. We also observe the initial inventory at the HF. Because local demand at the HF is not observed, we impute it as described below. We supplement these data with the estimated CHA-to-HF travel distances in miles:
\begin{equation}
    (d_1,d_2,\ldots,d_7)=(11, 3.8, 3, 0.25, 5, 5, 0.5).
    \label{eq:CHA_HF_distance_data}
\end{equation}

We focus on malaria cases among children under five. In Liberia, malaria cases are typically concentrated during the rainy season in spring and summer \citep{PMI_Liberia_Profile_2024}. In our data, virtually all malaria cases in this age group occur between first week in February and second week in July. Accordingly, we use this 22-week window as the planning horizon and focus on weekly review periods, that is, $T=22$ and $N=22$. 

We calibrate demand using the binomial demand model introduced in Section~\ref{sec:model}. The records of household visits imply an average of 28 households per CHA. For the HF local catchment, we use the Fangonda Public School voter-registration count as a population proxy \citep{NEC_Lofa_VRC_2023}. The voter-registration records report 438 registered voters, and assuming two registered voters per household gives 219 households. Under the assumption that each household has one child under five, we set the CHA and HF exposure sizes as $N_c=28$ and $N_h=219$, respectively. We argue below that our results are robust to this assumption.

\begin{figure}[!ht]
\centering
\includegraphics[width=0.7\textwidth]{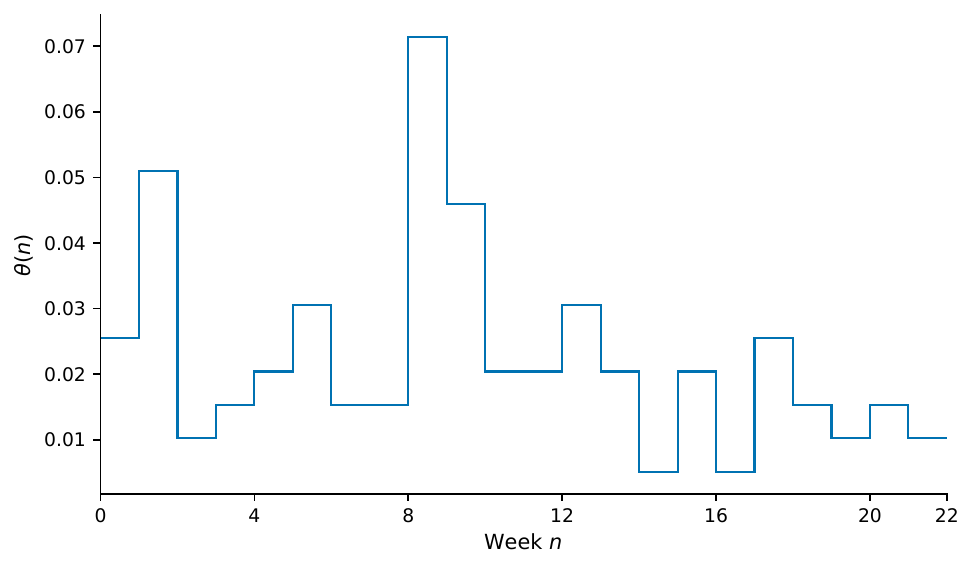}
\caption{Calibrated weekly incidence probability $\theta(n)$ for malaria demand among children under five in the Fangonda cluster.}
\label{fig:theta_calibration}
\end{figure}

Let $Y_k(n)$ denote the number of malaria cases identified by CHA $k$ among children under five in week $n$. We set the weekly incidence probability by pooling the CHA observations:
\begin{equation}
    \theta(n) = \frac{\sum_{k=1}^K Y_k(n)}{K N_c}, \quad n=1,\dots,N.
    \label{eq:weekly_incidence_probability}
\end{equation}
Figure~\ref{fig:theta_calibration} shows the resulting weekly incidence profile, which captures the seasonal pattern of malaria cases in the active period. Given $\theta(n)$, the weekly mean and variance inputs are
\begin{equation}
    \mu_k(n)=N_c\,\theta(n) \quad \text{and} \quad  \Gamma_{kk}(n)=N_c\,\theta(n)(1-\theta(n)), \quad k=1,\dots,K,
    \label{eq:weekly_demand_cha}
\end{equation}
and
\begin{equation}
    \mu_0(n) = N_h\,\theta(n) \quad \text{and} \quad \Gamma_{00}(n)=N_h\,\theta(n)(1-\theta(n)).
    \label{eq:weekly_demand_hf}
\end{equation}
In addition, we set off-diagonal covariance terms to zero initially, but consider dependence across locations in our sensitivity analysis. The estimate $\theta(n)$ in Equation~\eqref{eq:weekly_incidence_probability} relies on the assumption that each household has one child under the age of five on average. As shown in Figure~\ref{fig:theta_calibration}, $\theta(n)$ is a small number in the range of $0.01$--$0.07$. Now, consider an alternative assumption that, on average, a household has two children under the age of five. The new estimate, denoted by $\tilde{\theta}(n)$, is half of the earlier one, i.e., $\tilde{\theta}(n) = \theta(n) / 2$. Consider the corresponding mean and covariance estimates, denoted by $\tilde{\mu}$ and $\tilde{\Gamma}$, and note that $\tilde{\mu}_k(n) = 2 N_c \tilde{\theta}(n) = N_c \theta(n) = \mu_k(n)$ and
\begin{align*}
\tilde{\Gamma}_{kk}(n) &= 2N_c \tilde{\theta}(n) (1-\tilde{\theta}(n)) = N_c \theta(n) \Big(1-\frac{\theta(n)}{2}\Big) \approx N_c \theta(n)(1-\theta(n)) = \Gamma_{kk}(n),
\end{align*}
which follows because $\theta(n)$ is small.

We set the initial warehouse inventory to $Q(0)=142$, corresponding to the annual supply available to the Fangonda region in the data. Under the calibrated demand inputs, the total expected demand over the 22-week horizon is
\begin{equation}
    \bar{d} =  \sum_{n=1}^{N}\sum_{j=0}^{K}\mu_j(n) = 207.5,
\end{equation}
consisting of $109.5$ expected HF local cases and $98$ expected CHA cases. The available supply therefore covers only about $Q(0) / \bar{d} \approx 68.4\%$ of expected demand. Thus, the test problem represents a resource-constrained environment.

Next, we calibrate the stockout cost in Equation~\eqref{eq:cost_definition}. We first set the referral-compliance probability to $p_k=0.69$ for all CHAs \citep{brunner2022prereferral}\footnote{\cite{brunner2022prereferral} report that, in the Nigeria study site, $48\%$ of referred children under five with suspected severe malaria completed referral to a designated referral health facility, while an additional $21\%$ sought care from another public provider. Because our model contains a single referral-completion branch, we use the combined public-provider care-seeking rate as a baseline proxy for realized referral compliance.}. Meanwhile, we measure all costs in units of one untreated under-five malaria patient, so that $u=1$ in the objective. Let $w$ denote the caregiver's hourly opportunity cost and let $\kappa$ denote the one-way motorcycle-taxi fare per mile. For the untreated patients, we use evidence from a recent study of malaria episodes among children under five in Burundi, which reports that caregivers missed work for $8$ days over the duration of a malaria episode \citep{niyibitegeka2025economic}. Assuming an $8$-hour workday, the monetary proxy for one untreated patient is $8 \times 8 \times w$. For a referred patient from CHA $k$, we assume that visiting the HF requires one full day of caregiver time and a round-trip motorcycle-taxi ride. The referral cost is therefore
\[ r_k = \frac{8 w + 2 \kappa d_k}{64 w}, \]
where $d_k$ is the CHA-to-HF distance in Equation~\eqref{eq:CHA_HF_distance_data}. We set $w=\$0.68$ per hour as the minimum wage for workers in Liberia \citep{WageIndicatorLiberiaWorkWages2026}. We set $\kappa=\$0.15$ per mile, calculated from a reported fare of $\$2.50$ for approximately $27.4$ kilometers \citep{ShutterstockLiberiaMotorcycleTaxi2020}. 
Hence, the monetary value corresponding to one cost unit is $64w =43.52$, and the resulting normalized cost vector from Equation~\eqref{eq:cost_definition} is
\begin{equation}
    (c_0,c_1,\ldots,c_7) = (1,\, 0.4486,\, 0.4143,\, 0.4105,\, 0.3974,\, 0.4200,\, 0.4200,\, 0.3986).
    \label{eq:shortage_cost_calibration}
\end{equation}

We use $(\lambda,\alpha)=(0.25,0.5)$ as the reporting baseline, for which operational cost and weighted fairness utility have comparable magnitudes. In the additional experiments, we vary $\lambda$ while holding $\alpha=0.5$ fixed.

Finally, we choose the shipment bounds in each period using the calibrated demand mean and variance. For CHAs, we set
\begin{equation}
    \eta_k(n) = \chi\, \mu_k(n), \quad k=1,\ldots,K,\;n=1,2,\dots,N,
    \label{eq:eta_k}
\end{equation}
where $\chi>0$ scales the mean weekly load.
For the HF, the bound must account for both local HF demand and referrals from CHAs. We use $\beta = \max\{0,1 - Q(0) / \bar{d}\}$ as a crude estimate of the fraction of CHA demand that cannot be served locally, and therefore set
\begin{equation}
    \eta_0(n) = \chi\, \Big(\mu_0(n) + \beta \sum_{k=1}^K p_k \mu_k(n)\Big), \quad n=1,2,\dots,N.
    \label{eq:eta_0}
\end{equation}
As the baseline, we consider $\chi=1.25$.

Table \ref{tab:baseline_parameters} summarizes the model parameters used in the test problem.

\begin{table}[!ht]
\centering
\small
\caption{Model parameters in the test problem for the Fangonda cluster.}
\label{tab:baseline_parameters}
\begin{tabular}{lll}
\toprule
Parameter & Interpretation & Baseline Value \\
\midrule
$T$ & Planning horizon & $22$ \\
$K$ & Number of CHAs & $7$ \\
$\theta(n)$ & Weekly incidence probability & See Equation~\eqref{eq:weekly_incidence_probability} and Figure~\ref{fig:theta_calibration}\\
$\mu_k(n)$, $\Gamma_{kk}(n)$ & Weekly demand mean and variance & See Equations~\eqref{eq:weekly_demand_cha}--\eqref{eq:weekly_demand_hf}\\
$\Gamma_{k\ell}(n)$, $k\neq \ell$ & Cross-location covariance & $0$ \\
$p_k$ & Referral compliance probability & $0.69$ \\
$d_k$ & Distance from CHA $k$ to HF & $(11,3.8,3,0.25,5,5,0.5)$ \\
$u$ & Untreated cost unit & $1$ \\
$c_k$ & Shortage cost & See Equation~\eqref{eq:shortage_cost_calibration} \\
$Q(0)$ & Initial unallocated inventory & $142$ \\
$\lambda$ & Weight on fairness utility & $0.25$ \\
$\alpha$ & Curvature parameter of fairness utility & $0.5$ \\
$\eta_k$ & Shipment bound & See Equations~\eqref{eq:eta_k}--\eqref{eq:eta_0} \\
\bottomrule
\end{tabular}
\end{table}

\subsection{Benchmark Policies}
\label{sec:benchmarks}

We compare our policy with five benchmarks described as follows. The allocation decisions prescribed by these benchmarks will be adjusted as in Equation~\eqref{eq:hjb_discrete_feasible_allocation} if needed.

\paragraph{80/20 allocation.}
Under this policy, 80\% of the initial inventory is allocated to the HF and the remaining 20\% is divided equally among the CHAs. Subsequent shipments are then made to fulfill these fixed quotas at the maximum rate permitted per review period. Its shipment is 
\begin{align*}
    a_k(n)=\min\{\eta_k(n),(\bar A_k-A_k(n-1))^+\},\quad \bar A_0=0.8Q(0),\quad \bar A_k=0.2Q(0)/K,\ k=1,\dots,K.
\end{align*}
Interestingly, this 80/20 split corresponds to the stated allocation policy \citep{hrh2030_2021_liberia,kangar_2023_liberia,villagereach_2024_supplychain}, but our conversations with those running the system indicate that it is not often followed in practice. Our dataset also reflects this discrepancy. As we show below, the policy is ineffective.

\paragraph{Static policy.}
This policy follows the mean demand and ignores realized inventory states. Its shipment is $a_k(n)=\min\{\eta_k(n),\mu_k(n)\}$ for $k=0,1,\dots,K$. 

\paragraph{Fairness-only policy.}
This policy further accounts for the aggregate scarcity by scaling the static policy with a single factor
\[\gamma = \min\Bigg\{1,\frac{Q(0)}{\sum_{m=1}^{N}\sum_{\ell=0}^{K} \min\{\eta_\ell(m),\mu_\ell(m)\}}\Bigg\},\]
and the resulting shipment is 
\begin{equation}
    a_k(n)=  \gamma\min\{\eta_k(n),\mu_k(n)\}.
    \label{eq:scaled_static_policy}
\end{equation}
Under the current proportional shipment bounds, when the bounds do not bind, $\gamma=Q(0)/\bar d$ and every location, including the HF, receives the same fraction of its expected demand in each period. This spatially and temporally regular allocation maximizes the model's fairness utility under these conditions; see Appendix~\ref{app:deterministic_scaled_static_policy}.

\paragraph{Fluid-optimal policy.} This policy replaces demand by its mean and solves the resulting deterministic finite-horizon allocation problem. The solution gives an open-loop shipment schedule. The full formulation is given in Appendix~\ref{app:deterministic_fluid_policy}.

\paragraph{Myopic policy.}
This policy reacts only to current local inventory and ships to a location when its beginning-of-period inventory is at or below a threshold. Its shipment is 
\begin{align*}
    a_k(n)=\eta_k(n) \,\mathbf{1}\{Z_k(n-1)\le \epsilon\},\quad k=0,\dots,K.
\end{align*}
We set $\epsilon=0$ in our experiments.

The benchmark policies are implemented similarly to our proposed policy using Equation~\eqref{eq:hjb_discrete_feasible_allocation} to scale the allocations down if needed.

\section{Computational Results}
\label{sec:computational_results}

In this section, we first compare our policy with the benchmark policies introduced in the previous section. We then conduct additional experiments to examine the effects of varying (i) the weight $\lambda$ placed on fairness in the objective function, (ii) travel cost, (iii) spatial correlation in demand, (iv) supply scarcity, (v) HF demand, and (vi) demand variability.

As shown in Table~\ref{tab:baseline_compact_comparison}, our policy substantially outperforms the benchmarks overall. The strongest benchmark is the fairness-only policy, with a total objective value of $-9.736$. Relative to this benchmark, our policy reduces operational costs by about 6\% while achieving slightly lower fairness utility, approximately 98\% of that under the fairness-only policy. The static policy achieves slightly lower operational costs than our policy but yields substantially lower fairness utility. The 80–20 benchmark policy provides slightly higher overall coverage than our policy but much lower coverage for CHAs. It also performs worse than our policy in terms of both operational costs and fairness utility.

\begin{table}[!ht]
    \centering \small
    \setlength{\tabcolsep}{4pt}
    \renewcommand{\arraystretch}{1.12}
    \caption{Baseline policy comparison when $\lambda=0.25$.}
    \label{tab:baseline_compact_comparison}
    \begin{tabular}{lrrrrrr}
    \toprule
    Policy & \makecell{Total\\objective} & \makecell{Operating\\cost} & \makecell{Fairness\\utility} & Coverage & \makecell{Mean CHA\\local fill} & \makecell{HF shipment\\share} \\
    \midrule
    Our policy & -12.530 & 71.420 & 335.799 & 0.676 & 0.661 & 0.531 \\
    80--20 Rule & 12.300 & 72.799 & 241.994 & 0.688 & 0.293 & 0.800 \\
    Static & -0.909 & 70.683 & 286.370 & 0.679 & 0.672 & 0.528 \\
    Fairness-only & -9.736 & 76.090 & 343.307 & 0.655 & 0.624 & 0.528 \\
    Myopic & 19.827 & 77.857 & 232.118 & 0.652 & 0.495 & 0.635 \\
    Fluid-optimal & -2.744 & 81.899 & 338.569 & 0.618 & 0.781 & 0.365 \\
    \bottomrule
    \end{tabular}
\end{table}

\begin{figure}[!ht]
    \centering
    \includegraphics[width=\textwidth]{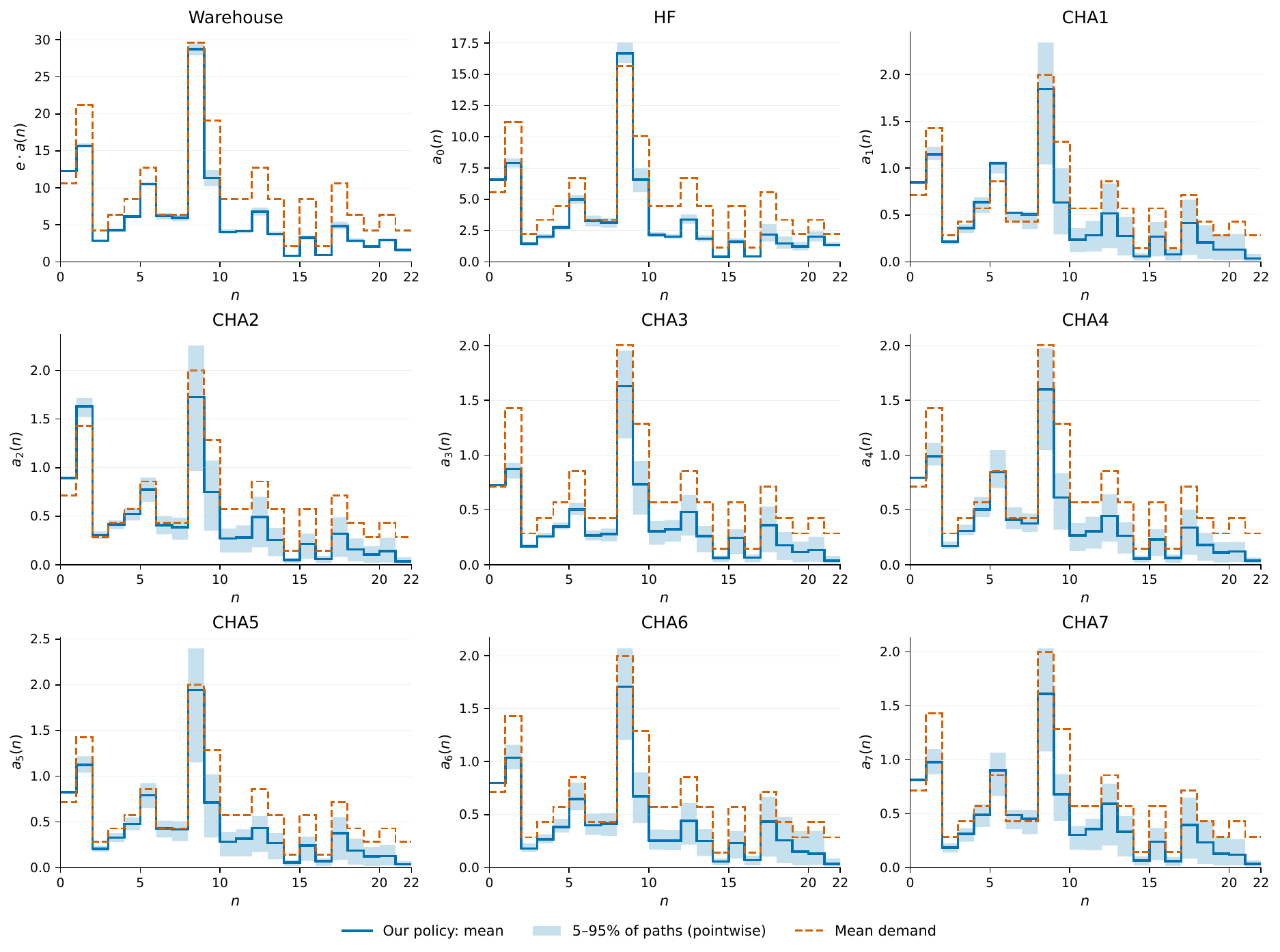}
    \caption{Weekly shipments under our policy in the baseline case. Blue lines show means over 10,000 demand paths, and shaded bands span the pointwise 5th--95th percentiles across paths. Orange dashed lines show mean weekly demand. The Warehouse panel compares total shipments with total local demand across all locations. Shipments and demand in week $n$ are displayed over $(n-1,n]$.}
    \label{fig:baseline_policy_trajectories}
\end{figure}

\begin{figure}[!ht]
    \centering
    \includegraphics[width=\textwidth]{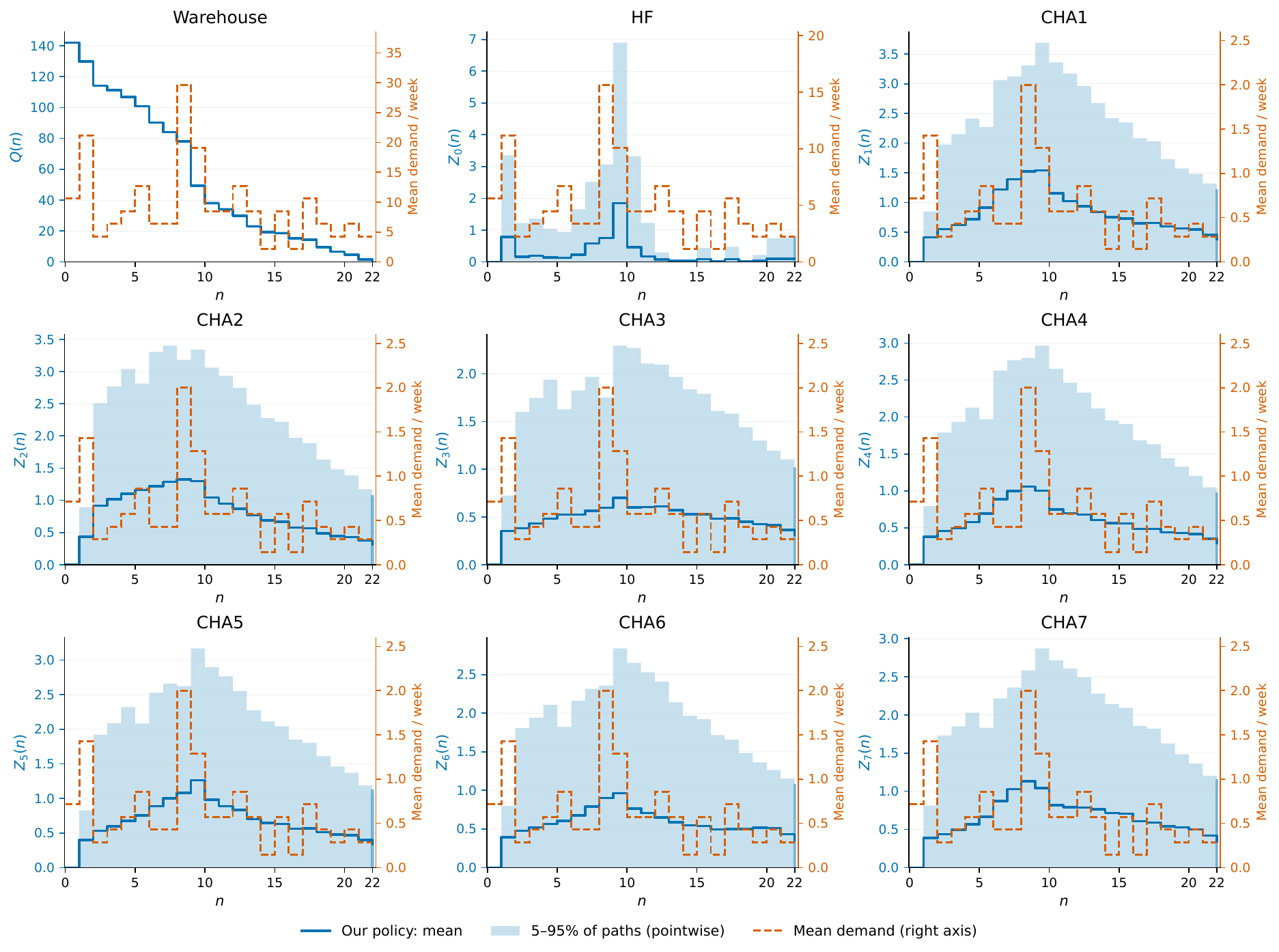}
    \caption{Inventory trajectories under our policy in the baseline case. Blue lines show means over 10,000 demand paths, and shaded bands span the pointwise 5th--95th percentiles across paths. Orange dashed lines show mean weekly demand on the right axis, using total local demand for Warehouse and excluding referrals for HF.}
    \label{fig:baseline_inventory_trajectories}
\end{figure}

Figures~\ref{fig:baseline_policy_trajectories} and~\ref{fig:baseline_inventory_trajectories} show the weekly shipments and inventory trajectories of our policy alongside mean demand in the baseline case. The shaded bands display variation across demand paths. Our policy releases 45.0\% of the initial warehouse inventory by the end of week 8, building local inventories ahead of the demand peak in week 9. After the peak, shipments generally decline and CHA inventories are gradually drawn down, leaving a mean total of 2.34 units at the CHAs at the end of week 22. The narrow warehouse inventory band contrasts with the wider bands at individual locations, indicating that a relatively stable aggregate release pattern coexists with substantial variation in local inventories across demand paths. Together with the performance comparison in Table~\ref{tab:baseline_compact_comparison}, these trajectories suggest that our policy improves on the fairness-only benchmark by building local inventory ahead of peak demand and adjusting subsequent shipments to leave less medication unused at the end of the planning horizon.

\FloatBarrier

\begin{figure}[!ht]
    \centering
    \includegraphics[width=0.6\textwidth]{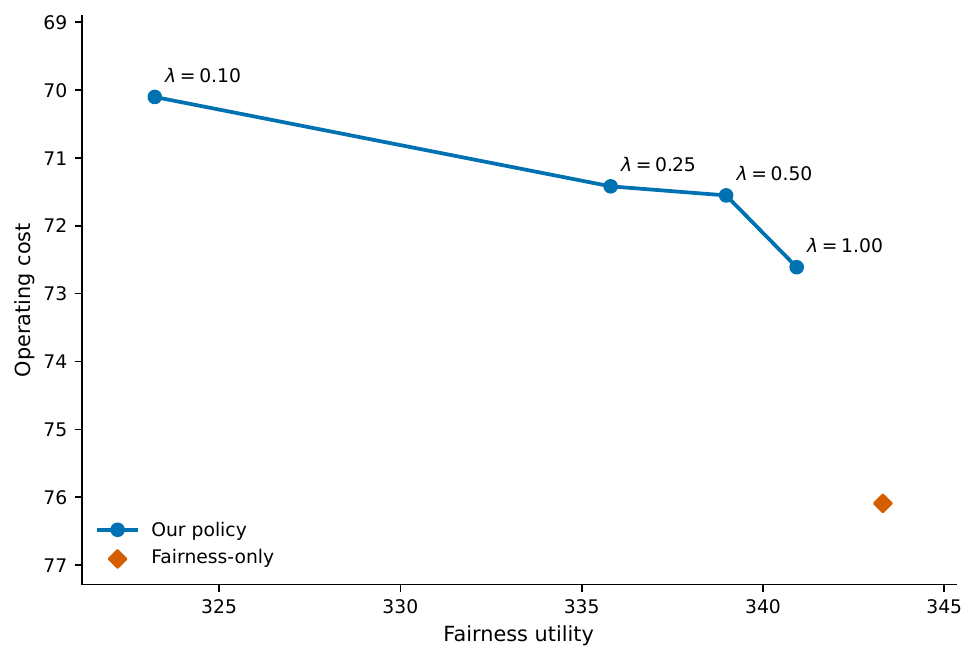}
    \caption{Cost--fairness trade-off over $\lambda\in\{0.10,0.25,0.50,1.00\}$ under the baseline calibration. The fairness-only policy is shown separately as the benchmark policy with maximum fairness.}
    \label{fig:fairness_weight_frontier}
\end{figure}

\paragraph{Varying the weight of fairness utility in the objective.} We consider $\lambda \in \{0.1, 0.25, 0.5, 1\}$. As $\lambda$ increases, our proposed policy places greater weight on allocation fairness. The results in Figure~\ref{fig:fairness_weight_frontier} suggest that $\lambda = 0.5$ may offer the best balance: increasing $\lambda$ from $0.25$ to $0.5$ yields a substantial increase in fairness utility with only a slight increase in operational costs.

\paragraph{Varying the travel cost.} We vary the one-way transportation fare per mile, $\kappa \in \{0.05, 0.15, 0.25\}$. Higher values of $\kappa$ increase the cost of completed referrals and make local treatment at CHAs more valuable. Table~\ref{tab:travel_cost_summary} presents the results.

\begin{table}[!ht]
    \centering \small
    \setlength{\tabcolsep}{4pt}
    \renewcommand{\arraystretch}{1.12}
    \caption{Travel-cost response under $\lambda=0.25$.}
    \label{tab:travel_cost_summary}
    \begin{tabular}{rrrrrrr}
    \toprule
    $\kappa$ & \makecell{CHA shipment\\share} & \makecell{CHA local\\unmet} & \makecell{Release fraction\\through week 8} & \makecell{Mismatch\\frequency} & \makecell{Mean absolute\\severity} \\
    \midrule
    0.05 & 0.405 & 42.160 & 0.510 & 0.884 & 1.399 \\
    0.15 & 0.469 & 33.720 & 0.450 & 0.852 & 1.276 \\
    0.25 & 0.531 & 25.090 & 0.430 & 0.779 & 1.033 \\
    \bottomrule
    \end{tabular}
\end{table}

As expected, higher travel costs lead to larger inventory allocations to CHAs, reducing both the amount of CHA demand that is not met locally and the number of referrals. A mismatch event occurs when one CHA runs out of medication while another has excess inventory in the same period. Both the frequency and magnitude of these local supply–demand mismatches decline as travel costs increase. However, the overall rate at which inventory is allocated over time does not increase, because higher allocations to CHAs are accompanied by lower allocations to the HF.

\paragraph{Increasing Demand Variability.} To assess the robustness of our findings to demand variability, we increase the variance of demand at each CHA while holding its mean fixed. Letting $\tilde{\Gamma}$ denote the new diagonal covariance matrix, we set
\[ \tilde{\Gamma}_{kk}(n) = c\,\Gamma_{kk}(n), \qquad k = 1,\ldots,K, \]
where $c \in \{1,2,5\}$. The case $c=1$ corresponds to our original test problem. For $c \in \{2,5\}$, we assume that demand at CHA $k$ follows a beta-binomial distribution with variance $\tilde{\Gamma}_{kk}(n)$ and the original mean $\mu_k(n)$; see Appendix~\ref{app:demand_generation} for further details.

\begin{table}[!ht]
    \centering
    \small
    \setlength{\tabcolsep}{4pt}
    \renewcommand{\arraystretch}{1.12}
    \caption{Demand-overdispersion response when $\lambda=0.25$.}
    \label{tab:demand_overdispersion_performance}
    \begin{tabular}{lrrrrrrr}
    \toprule
    Variance & \makecell{Objective gap vs.\\fairness-only} & \makecell{Operating\\cost} & \makecell{Fairness\\utility} & Coverage & \makecell{Mismatch\\frequency} & $\bar S$ & $S_+$ \\
    \midrule
    1$\times$ & -2.794 & 71.420 & 335.799 & 0.676 & 0.852 & 1.276 & 1.497 \\
    2$\times$ & -4.698 & 73.660 & 334.153 & 0.668 & 0.725 & 1.458 & 2.012 \\
    5$\times$ & -7.470 & 79.714 & 332.506 & 0.647 & 0.536 & 1.861 & 3.474 \\
    \bottomrule
    \end{tabular}
    \par\vspace{2pt}
    \parbox{0.98\textwidth}{\footnotesize \textit{Note:} $\bar S$ is mean absolute mismatch severity and $S_+$ is mean severity conditional on a positive mismatch event.}
\end{table}
    
Table~\ref{tab:demand_overdispersion_performance} summarizes the results. As demand variance increases, operating costs rise and demand coverage declines. However, our policy's improvement over the strongest benchmark, the fairness-only policy, also increases. Although mismatches occur less frequently, their magnitudes increase substantially.

\paragraph{Increasing Spatial Demand Correlations.} We consider $\rho\in\{0, 0.5, 0.7\}$ for the pairwise correlation in weekly demand across the HF and CHAs. To model these correlations, we use a calibrated Gaussian copula with binomial marginal distributions, preserving the mean and variance of demand at each location; see Appendix~\ref{app:demand_generation} for details.

\begin{table}[!ht]
    \centering
    \small
    \setlength{\tabcolsep}{4pt}
    \renewcommand{\arraystretch}{1.12}
    \caption{Demand-correlation response when $\lambda=0.25$.}
    \label{tab:demand_correlation_summary}
    \begin{tabular}{lrrrrrrr}
    \toprule
    $\rho$ & \makecell{Objective gap vs.\\fairness-only} & \makecell{Operating\\cost} & \makecell{Fairness\\utility} & \makecell{Mean CHA\\local fill} & \makecell{CHA shipment\\share} & \makecell{Mismatch\\frequency} & \makecell{Mean absolute\\severity} \\
    \midrule
    0.00 & -2.794 & 71.420 & 335.799 & 0.661 & 0.469 & 0.852 & 1.276 \\
    0.50 & -2.395 & 71.080 & 330.293 & 0.697 & 0.480 & 0.547 & 0.591 \\
    0.70 & -2.868 & 71.338 & 332.842 & 0.757 & 0.526 & 0.405 & 0.384 \\
    \bottomrule
    \end{tabular}
\end{table}

Table~\ref{tab:demand_correlation_summary} shows that our policy outperforms the benchmarks across all values of $\rho$ considered. The CHAs' share of inventory allocations increases with $\rho$. Recall that a mismatch event occurs when one location experiences a shortage while another has excess inventory in the same period. As expected, both the frequency and magnitude of these mismatches decline as $\rho$ increases.

\paragraph{Varying HF Demand.} Recall that local HF demand is unobserved, so we impute it using Fangonda voter-registration data, assuming two registered voters per household. We vary this assumption by considering $N_h \in \{199,219,243\}$, corresponding respectively to $2.2$, $2$, and $1.8$ registered voters per household. For each value of $N_h$, we recompute the mean HF demand $\mu_0(n)$, its variance $\Gamma_{00}(n)$, and the HF shipment bound $\bar{a}_0(n)$. This experiment changes both total expected demand and its composition between the HF and CHAs, while leaving CHA demand and initial inventory fixed.

\begin{table}[!ht]
    \centering
    \small
    \setlength{\tabcolsep}{4pt}
    \renewcommand{\arraystretch}{1.12}
    \caption{HF-demand response when $\lambda=0.25$.}
    \label{tab:hf_demand_summary}
    \begin{tabular}{rrrrrrr}
    \toprule
    $N_h$ & \makecell{Total expected\\demand} & \makecell{Objective gap vs.\\fairness-only} & \makecell{HF shipment\\share} & Coverage & \makecell{HF-load\\fulfillment} & \makecell{Mean CHA\\local fill} \\
    \midrule
    199 & 197.5 & -3.724 & 0.500 & 0.709 & 0.597 & 0.706 \\
    219 & 207.5 & -2.794 & 0.531 & 0.676 & 0.572 & 0.661 \\
    243 & 219.5 & -0.827 & 0.550 & 0.644 & 0.539 & 0.645 \\
    \bottomrule
    \end{tabular}
\end{table}

Table~\ref{tab:hf_demand_summary} shows that the policy shifts inventory toward the HF as HF demand increases. Because total expected demand also rises while initial inventory remains fixed, coverage and both facility-specific service measures decline. Our policy continues to be the best policy though its improvement over the fairness-only policy decreases as HF demand increases.

\paragraph{Varying Initial Inventory.} We vary the initial inventory $Q(0) \in \{124.5,142,166\}$, corresponding to $60.0\%$, $68.4\%$, and $80.0\%$ of total expected demand, while holding all other parameters fixed. The middle value, $Q(0)=142$, is the base case. Table~\ref{tab:initial_inventory_response} presents the results.

\begin{table}[!ht]
    \centering
    \small
    \setlength{\tabcolsep}{4pt}
    \renewcommand{\arraystretch}{1.12}
    \caption{Initial-inventory response when $\lambda=0.25$.}
    \label{tab:initial_inventory_response}
    \begin{tabular}{rrrrrrrr}
    \toprule
    $Q(0)$ & \makecell{Objective gap vs.\\fairness-only} & \makecell{Week-8 release,\\ units (\%)} &  Coverage & \makecell{HF-load\\fulfillment} & \makecell{HF shipment\\share} & \makecell{HF\\shipment} & \makecell{CHA\\shipment}\\
    \midrule
    124.5 & 1.873 & 63.851 (51.3\%) & 0.598 & 0.451 & 0.477 & 59.331 & 65.169\\
    142   & -2.794 & 63.884 (45.0\%)  & 0.676 & 0.572 & 0.531 & 75.441 & 66.557 \\
    166   & -4.653 & 73.706 (44.4\%) & 0.790 & 0.745 & 0.587 & 97.491 & 68.491 \\
    \bottomrule
    \end{tabular}
\end{table}

\begin{figure}[!ht]
    \centering
    \includegraphics[width=0.49\textwidth]{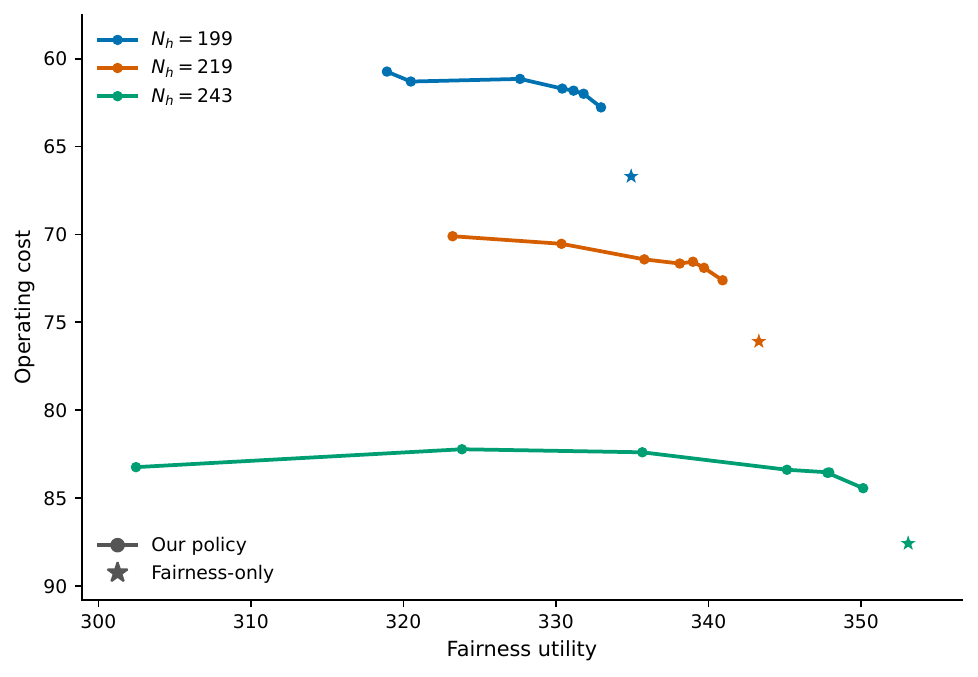}\hfill
    \includegraphics[width=0.49\textwidth]{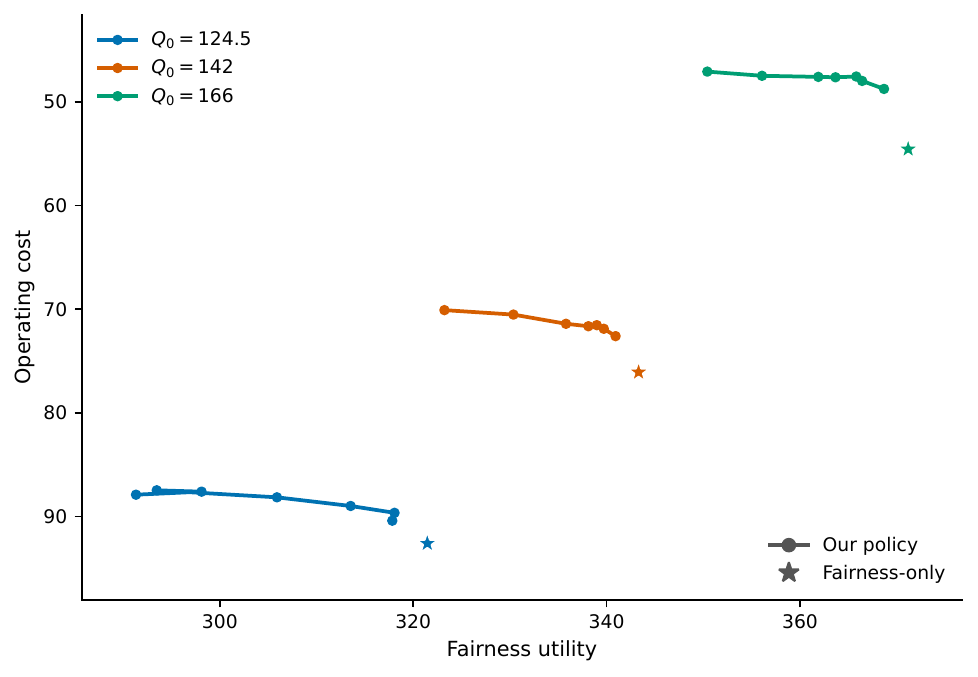}
    \caption{Operating-cost and raw-fairness trade-offs across fairness weights and scenarios. Circles show our policy for $\lambda=0.10,0.15,0.25,0.35,0.50,0.70,1.00$, and stars show fairness-only policy. The left panel varies HF-demand exposure with $Q(0)=142$, and the right panel varies initial inventory with $N_h=219$.}
    \label{fig:common_raw_lambda_frontier}
\end{figure}

Table~\ref{tab:initial_inventory_response} shows that increasing initial inventory substantially improves service performance, and changes the relative performance of our policy against the fairness-only policy from unfavorable under low inventory to increasingly favorable at higher inventory levels. This follows because when the initial inventory is low (or very high) there is little room to improve the operating cost. Figure~\ref{fig:common_raw_lambda_frontier} displays the operating cost--fairness utility tradeoff for different levels of HF demand (left panel) and for different levels of initial inventory (right panel). As the fairness weight $\lambda$ increases, the proposed policy results in higher fairness utility, approaching that under the fairness-only policy---albeit incurring a higher operating cost.

In addition, when the initial inventory increases, the additional inventory is allocated mainly to the HF. Although shipments to both the HF and the CHAs increase with $Q(0)$, the increase in HF shipments is much larger, resulting in a higher HF shipment share. 
Interestingly, lower initial inventory does not reduce the total quantity of inventory allocated to the CHAs and the HF early in the planning horizon. The low-inventory and base cases release nearly the same amount by week 8, so the low-inventory policy commits a larger fraction of its available stock early and depletes it sooner.

\section{Conclusion}
\label{sec:conclusion}
We study the dynamic allocation of scarce medications between a health facility and the Community Health Assistants it supervises in Liberia's community health program. The allocation problem requires balancing the flexibility of retaining inventory centrally against the benefit of placing medications closer to patients, while also accounting for geographical fairness considerations. We formulate this problem as a finite-horizon, two-echelon stochastic control problem that captures uncertain demand, referrals, and the costs associated with untreated patients and referrals.

To solve this problem, we develop a continuous-time approximation based on reflected Brownian motion and characterize the allocation decision through a high-dimensional HJB equation. We then adapt the deep splitting method to approximate the value function and the marginal values that determine the allocation policy. This approach allows us to solve the resulting high-dimensional dynamic allocation problem without discretizing the full state space.

We calibrate the model using data from the Fangonda region of Lofa County, Liberia, and evaluate the proposed policy in the original discrete-review system. Our policy outperforms benchmark policies considered in the paper. Relative to the fairness-only policy, the strongest benchmark, our policy reduces operating costs by about 6\% while retaining approximately 98\% of its fairness utility. These results show that dynamically adjusting allocations based on the system state can improve the balance between operating costs and allocation fairness significantly.

Our additional experiments show how the proposed policy responds to changes in the operating environment. Increasing the weight placed on fairness moves the policy toward higher fairness utility at the cost of somewhat higher operating costs. Higher referral costs shift inventory toward CHAs and reduce local unmet demand and supply--demand mismatches. Greater demand variability makes the system more difficult to serve, but also increases the improvement of our policy over the fairness-only benchmark. Higher spatial correlation reduces mismatches across locations, while changes in HF demand and initial inventory lead the policy to adjust the balance of shipments between the HF and CHAs. Together, these experiments show that the proposed policy responds systematically to changes in costs, demand, and resource availability.

Overall, our results suggest that dynamic allocation can improve the use of scarce medication inventory by accounting for current inventory states, uncertain demand, referral costs, and fairness considerations. More broadly, the model and computational approach provide a way to study similar last-mile allocation problems in which limited resources must be distributed dynamically across multiple service locations.

\noindent\textbf{Code and data availability.}
Code and computational materials for reproducing the figures and tables
are available at \url{https://github.com/YuwanIvan/LiberiaRBM_public}.

\noindent\textbf{Acknowledgements.} The authors thank the staff at the International Rescue Committee for facilitating the discussions and data compilation essential to this project, with special thanks to Dianah Majekodunmi, Rukmini Giriharadas, and Erika Pearl. We also thank participants at the 2026 INFORMS Healthcare Conference for their helpful comments and suggestions. We would also like to acknowledge the International Rescue Committee's support of J.R. Montgomery, who worked on this research while serving as a postdoctoral research fellow with the organization. Any opinions, findings, and conclusions or recommendations expressed in this material are those of the authors and do not necessarily reflect the views of the International Rescue Committee.

\bibliography{ref}
\clearpage
\setstretch{1.45}

\appendix

\section{Details of Benchmark Policies}

\subsection{Fairness-only Policy}
\label{app:deterministic_scaled_static_policy}
Consider the fairness-only problem that ignores stockout costs and state dynamics, with nonbinding shipment bounds:
\begin{align*}
\max_{a_k(n) \geq 0} \quad  \sum_{n=1}^N \Psi_\alpha(f(n);\mu(n)) \qquad  \text{s.t.} \quad  \sum_{n=1}^N\sum_{k=0}^K a_k(n) \le Q(0).
\end{align*}
Using Equation~\eqref{eq:fairness_utility} and H\"{o}lder's inequality, we have
\begin{align*}
(1-\alpha) \sum_{n=1}^N \Psi_{\alpha}(f(n);\mu(n)) \leq \Big(\sum_{n=1}^N \sum_{k=0}^K a_k(n) \Big)^{1-\alpha} \Big(\sum_{n=1}^N \sum_{k=0}^K \mu_k(n)\Big)^{\alpha} \leq [Q(0)]^{1-\alpha} \bar{d}^{\alpha},
\end{align*}
where $\bar{d} = \sum_{n=1}^N \sum_{k=0}^K \mu_k(n)$ is the mean total demand, and the equality is attained when
\begin{equation}
a_k(n) = \frac{Q(0)}{\bar{d}} \mu_k(n).
\label{eq:max_fairness}
\end{equation}
Hence, when $Q(0)<\bar d$ and $\eta_k(n)\geq\mu_k(n)$ for every location-period pair, Equation~\eqref{eq:scaled_static_policy} attains this maximum. Each location receives the fraction $Q(0)/\bar d$ of its expected demand in every period. This motivates the name of the policy.

\subsection{Fluid-optimal Policy}
\label{app:deterministic_fluid_policy}
The fluid-optimal policy replaces random demand by its calibrated mean and solves a deterministic finite-horizon allocation problem. The control variables are the allocations $a_k(n)$. We introduce auxiliary variables $Q(n)$, $Z_k(n)$, and $\ell_k(n)$ to express the warehouse balance, local inventory balance, and unmet demand terms linearly. Given $\alpha = 0.5$, the resulting optimization problem is
\begin{equation}
\begin{aligned}
\min_{a,Q,Z,\ell} \quad & \sum_{n=1}^N \sum_{k=0}^K c_k \ell_k(n) - 2 \lambda \sum_{n=1}^N \sum_{k=0}^K \sqrt{a_k(n)\mu_k(n)} \\
\text{s.t.} \quad & Q(n) = Q(n-1) - \sum_{k=0}^K a_k(n),\quad n=1,\ldots,N, \\
& Z_k(n) = Z_k(n-1) + a_k(n) - \mu_k(n) + \ell_k(n), \quad n=1,\ldots,N,\ k=1, \ldots,K, \\
& Z_0(n) = Z_0(n-1) + a_0(n) - \mu_0(n) - \sum_{k=1}^K p_k\ell_k(n) + \ell_0(n), \quad n=1,\ldots,N, \\
& a_k(n),Q(n),Z_k(n),\ell_k(n)\ge 0, \quad n=1,\ldots,N,\ k=0,\ldots,K, \\
& a_k(n)\le \eta_k(n) h, \quad n=1,\ldots,N,\ k=0,\ldots,K, \\
& Q(0)=q_0,\quad Z_k(0)=0, \ k=0,\ldots,K.
\end{aligned}
\label{eq:fluid_optimal_problem_general}
\end{equation}
When $\lambda=0$, Equation~\eqref{eq:fluid_optimal_problem_general} is a linear program (LP). When $\lambda>0$, we further introduce auxiliary variables $y_k(n)$ such that
\begin{equation*}
    y_{k}(n)\geq 0,\quad [y_k(n)]^2 \leq a_k(n)\mu_k(n).
\end{equation*} 
Then, Equation~\eqref{eq:fluid_optimal_problem_general} is equivalent to a constrained quadratic program (CQP): 
\begin{align*}
    \min_{a,Q,Z,\ell,y} \quad & \sum_{n=1}^N \sum_{k=0}^K c_k \ell_k(n) - 2 \lambda \sum_{n=1}^N \sum_{k=0}^K y_k(n) \\
    \text{s.t.} \quad & (a,Q,Z,\ell) \in \mathcal{C}(q_0)\\
    &  [y_k(n)]^2 \leq a_k(n) \mu_k(n), \quad n=1,\ldots,N,\ k=0,\ldots,K,\\
    & y_k(n) \ge 0, \quad n=1,\ldots,N,\ k=0,\ldots,K,
\end{align*}
where $\mathcal{C}(q_0)$ is the feasible set for $(a,Q,Z,\ell)$ in Equation~\eqref{eq:fluid_optimal_problem_general}. In the implementation, we solve both the LP case and the CQP case using Gurobi.

\section{Evaluation Metrics}
\label{app:evaluation_metrics}

In this appendix, we define the evaluation metrics used in Section~\ref{sec:computational_results}. Unless stated otherwise, we calculate each metric on a simulated path over review periods $n=1,\ldots,N$ and then average across paths. The HF has index $0$ and the CHAs have indices $1,\ldots,K$; the path index is suppressed below except where needed.

For each setting, we evaluate our policy and benchmark policies on $R=10{,}000$ demand paths. If $X_r$ is a pathwise metric, its Monte Carlo mean and standard error are
\[ \hat X=\frac{1}{R}\sum_{r=1}^R X_r,\quad \operatorname{se}(\hat X)=\left[\frac{1}{R(R-1)}\sum_{r=1}^R(X_r-\hat X)^2\right]^{1/2}. \]

Let $\Delta L_k(n)=L_k(n)-L_k(n-1)$. For $k\geq1$, $\Delta L_k(n)$ is \emph{CHA local unmet} demand before possible referral; $\Delta L_0(n)$ is unmet demand at the HF after referral arrives. Our model adds the probability-weighted referral demand $p_k\Delta L_k(n)$ from each CHA to the HF, giving HF load
\[ H(n)=\xi_0(n)+\sum_{k=1}^K p_k\Delta L_k(n). \]
Let $\operatorname{Serv}_k(n)$ denote service at location $k$ in week $n$. HF service includes treatment of referred patients, whereas CHA service is local:
\[ \operatorname{Serv}_0(n)=H(n)-\Delta L_0(n),\quad \operatorname{Serv}_k(n)=\xi_k(n)-\Delta L_k(n),\; k=1,\ldots,K. \]

\paragraph{Objective.}
For evaluation, we report the pathwise objective underlying \eqref{eq:objective_discrete}, denoted by $J$. In addition, we decompose $J$ into cost and fairness as
\begin{align*}
    J = C - \lambda U \quad \text{with} \quad  C = c\cdot L(N) \quad \text{and} \quad U = \sum_{n=1}^{N} \Psi_{\alpha}(f(n);\mu(n)),
\end{align*}
where $C$ is the operational cost, and $U$ is the fairness utility. Furthermore, we decompose the operational cost $C$ into untreated cost and referral cost defined as $C = C_{\mathrm u} + C_{\mathrm{ref}}$, where
\[ C_{\mathrm u}=u\sum_{n=1}^N\left[\Delta L_0(n)+\sum_{k=1}^K(1-p_k)\Delta L_k(n)\right], \quad C_{\mathrm{ref}}=\sum_{n=1}^N\sum_{k=1}^K p_k r_k\Delta L_k(n). \]
We also compare the objective of our policy $J$ with that of fairness-only policy $J_r^{\mathrm{FO}}$, and we report the objective gap between them defined as $\Delta J = J - J^{\mathrm{FO}}$.

\paragraph{Service.}

We report three coverage measures. A ratio with a zero pathwise demand denominator is set to zero before averaging. Overall coverage is the fraction of demand treated at a CHA or the HF, including treatment after referral:
\[ D=\sum_{n=1}^N\sum_{k=0}^K\xi_k(n),\qquad
\mathrm{Cov}=\frac{\sum_n\sum_{k=0}^K\operatorname{Serv}_k(n)}{D}. \]
CHA coverage is the fraction of CHA-origin demand treated locally, before referral. With $D_{\mathrm{CHA}}=\sum_{n,k\ge1}\xi_k(n)$, it is
\[ \mathrm{LF}=\frac{\sum_{n,k\ge1}\operatorname{Serv}_k(n)}{D_{\mathrm{CHA}}}. \]
This ratio pools demand across CHAs; patients treated at the HF after referral do not count toward CHA coverage. HF coverage is the fraction of HF load treated there, including local HF demand and referral demand in the denominator:
\[ \mathrm{HFCov}=\frac{\sum_n \operatorname{Serv}_0(n)}{\sum_n H(n)}. \]

\paragraph{Allocation.} We describe where and when the warehouse releases inventory. Total shipment and its destination shares are
\[ e\cdot A(N)=\sum_{n,k}a_k(n),\qquad
\mathrm{HFShare}=\frac{A_0(N)}{e\cdot A(N)},\qquad
\mathrm{CHAShare}=\frac{\sum_{k=1}^K A_k(N)}{e\cdot A(N)}. \]
The two shares describe the HF--CHA split of executed shipments. To describe timing, cumulative release through week $m$ is $e\cdot A(m)=\sum_{n=1}^m\sum_{k=0}^K a_k(n)$, and its fraction of initial warehouse inventory is $e\cdot A(m)/Q(0)$. We report this fraction through week 8; it does not indicate where the released inventory went.

\paragraph{Inventory use and spatial mismatch.}
End-of-horizon utilization is the fraction of initial system stock consumed:
\[ \mathrm{Util}=1-\frac{Q(N)+\sum_{k=0}^K Z_k(N)}{Q(0)+\sum_{k=0}^K Z_k(0)}. \]
To identify stock left at one CHA while demand goes unmet at another, consider week $n$ on path $r$. Let $I_{n,r}=\sum_{k=1}^K Z_k^{(r)}(n)$ be remaining CHA inventory and $B_{n,r}=\sum_{k=1}^K\Delta L_k^{(r)}(n)$ be CHA local unmet demand. Define mismatch severity $S_{n,r}=\min\{I_{n,r},B_{n,r}\}$ and event indicator $M_{n,r}=\mathbf 1\{S_{n,r}>0\}$. For event detection, stock and local unmet demand count as positive only above $10^{-12}$. Because a CHA with unmet local demand has no stock left, a positive $S_{n,r}$ reflects stock and unmet demand at different CHAs.

We report mismatch frequency and mean severity across all path-weeks, plus mean severity conditional on a mismatch:
\[\mathrm{MismatchFreq}=\frac{1}{RN}\sum_{r=1}^R\sum_{n=1}^N M_{n,r},\qquad
\bar S=\frac{1}{RN}\sum_{r=1}^R\sum_{n=1}^N S_{n,r},\]
\[S_+=\frac{\sum_{r=1}^R\sum_{n=1}^N S_{n,r}M_{n,r}}{\sum_{r=1}^R\sum_{n=1}^N M_{n,r}}.\]
The $S_+$ denominator pools positive events across paths and weeks; it is undefined if none occurs. Mismatch severity represents hypothetical same-week pairing of CHA stock and unmet demand, not feasible transfers or avoidable untreated demand.

\section{Derivation of the Deep Splitting Algorithm}
\label{appendix:derivation_of_deep_splitting_algorithm}
In this appendix, we provide a formal derivation of the deep splitting algorithm in Section \ref{sec:computational_method}. We freeze the nonlinear term on each fine interval and condition on the reference state to obtain the regression target.

To derive the algorithm, we focus on the following equivalent characterization of the HJB equation in \eqref{eq:hjb_finite} for any $t\in[0,T]$, $z\in\R^{K+1}$ and $q\in\R_{+}$:
\begin{equation}
\begin{aligned}
\frac{\partial}{\partial t} V(t,z,q) + \mathcal{L}_t V(t,z,q) + [\tilde{\varphi}\cdot \mathcal{D} V] (t,z,q) + F(t, z, q, \mathcal{D} V(t, z, q)) = 0,
\end{aligned}
\label{eq:hjb_finite_F}
\end{equation}
where $F$ is defined in \eqref{eq:F_def}.

First, we fix a partition of $[0,T]$ such that $0 = \tau_0 < \tau_1 < \dots < \tau_M = T$ with $\Delta \tau_m =\tau_{m+1} - \tau_{m}$ for $m=0,1,\dots,M-1$. Under suitable integrability assumptions, for $m=0,1,\dots,M-1$, $t \in [\tau_m, \tau_{m+1}]$, $z \in \mathbb{R}^{K+1}$, and $q \in \mathbb{R}_{+}$, Equation \eqref{eq:hjb_finite_F} implies that
\begin{equation}
    \begin{aligned}
         V(t, z, q) &= V(\tau_m, z, q) + \int_{\tau_m}^{t} \frac{\partial V}{\partial t}(s, z, q) \dif s\\
        & = V(\tau_m, z, q) - \int_{\tau_m}^{t} [\mathcal{L}_s V + \tilde{\varphi}\cdot \mathcal{D} V] (s,z,q) \dif s - \int_{\tau_m}^{t} F(s, z, q, \mathcal{D} V(s, z, q)) \dif s.
    \end{aligned}
\end{equation}
For a sufficiently fine partition, one can write
\begin{equation}
    \begin{aligned}
        & V(t, z, q)  \approx V(\tau_m, z, q) - \int_{\tau_m}^{t}  [\mathcal{L}_s V + \tilde{\varphi}\cdot \mathcal{D} V] (s,z,q)  \dif s -  F(\tau_m, z, q, \mathcal{D} V(\tau_m, z, q))  \Delta \tau_m.
    \end{aligned}
    \label{eq:V_linear_approx}
\end{equation}

We assume there exists a function $\mathcal{V}: [0, T] \times \mathbb{R}^{K+1}_{+} \times \mathbb{R}_{+} \to \mathbb{R}$ satisfying $\mathcal{V}(T, z, q) = 0$ for any $z \in \mathbb{R}^{K+1}_{+}$ and $q \in \mathbb{R}_{+}$, such that for every $m = 0, 1, \dots, M-1$, $\mathcal{V}|_{(\tau_m,\tau_{m+1}] \times \mathbb{R}^{K+1}_{+} \times \mathbb{R}_{+}}$ belongs to $C^{1,2,1}((\tau_m,\tau_{m+1}] \times \mathbb{R}^{K+1}_{+} \times \mathbb{R}_{+})$ with at most polynomially growing partial derivatives. In addition,  $\mathcal{V}(t, z, q)$ satisfies
\begin{equation}
\begin{aligned}
    \mathcal{V}(t, z, q) & \approx \mathcal{V}(\tau_m, z, q) - \int_{\tau_m}^{t} [\mathcal{L}_s \mathcal{V} + \tilde{\varphi}\cdot \mathcal{D} \mathcal{V}] (s,z,q)  \dif s - F(\tau_m, z, q, \mathcal{D} \mathcal{V}(\tau_m, z, q))  \Delta \tau_m.
\end{aligned}
\label{eq:V_linear_approx_eq}
\end{equation}
for all $m = 0, 1, \dots, M-1$, $t \in (\tau_m,\tau_{m+1}]$, $z \in \mathbb{R}^{K+1}_{+}$, and $q \in \mathbb{R}_{+}$. The preceding equation is essentially the aforementioned linear approximation of the HJB equation \eqref{eq:hjb_finite_F} over small time intervals $(\tau_m,\tau_{m+1}]$ for $m = 0, 1, \dots, M-1$. Comparing \eqref{eq:V_linear_approx_eq} with \eqref{eq:V_linear_approx}, one expects that
\begin{equation}
    \begin{aligned}
        \mathcal{V}(t, z, q) &\approx V(t, z, q).
    \end{aligned}
    \label{eq:V_approx}
\end{equation}
It follows from \eqref{eq:V_linear_approx_eq} that, for $m = 0, 1, \dots, M-1$,
\begin{equation}
    \begin{aligned}
        \frac{\partial}{\partial t} \mathcal{V}(t,z,q) + \mathcal{L}_t \mathcal{V}(t,z,q) + [\tilde{\varphi}\cdot \mathcal{D} \mathcal{V}] (t,z,q) = 0, \quad \forall (t, z, q) \in (\tau_m,\tau_{m+1}]\times \mathbb{R}^{K+1}\times \mathbb{R}_{+},
    \end{aligned}
    \label{eq:hjb_finite_F_linear}
\end{equation}
From It\^{o}'s formula for RBM \citep{harrison1981reflected}, we have
\begin{equation}
    \begin{aligned}
        \dif \mathcal{V}(t, \tilde{Z}(t),\tilde{Q}(t)) &= \Big(\frac{\partial}{\partial t} + \mathcal{L}_t \Big)  \mathcal{V}(t,\tilde{Z}(t),\tilde{Q}(t)) \dif t + [\tilde{\varphi} \cdot \mathcal{D} \mathcal{V}](t,\tilde{Z}(t),\tilde{Q}(t)) \dif t\\
        & \qquad \qquad + \nabla_z \mathcal{V}(t,\tilde{Z}(t),\tilde{Q}(t)) \cdot \dif W(t) + \mathcal{R} \mathcal{V}(t,\tilde{Z}(t),\tilde{Q}(t)) \cdot \dif \tilde{L}(t).
    \end{aligned}
    \label{eq:ito_formula_linear}
\end{equation}
Substituting \eqref{eq:hjb_finite_F_linear} and the boundary condition \eqref{eq:hjb_boundary_condition_finite} into \eqref{eq:ito_formula_linear} yields
\begin{equation}
\begin{aligned}
\dif \mathcal{V}(t, \tilde{Z}(t),\tilde{Q}(t)) &= \nabla_z \mathcal{V}(t,\tilde{Z}(t),\tilde{Q}(t)) \cdot \dif W(t) - c \cdot \dif \tilde{L}(t).
\end{aligned}
\label{eq:equivalent_sde_finite_linear}
\end{equation}
Then integrating both sides on $(t, \tau_{m+1}]$ for $t \in (\tau_m, \tau_{m+1}]$ gives that $\mathbb{P}$-a.s.
\begin{equation}
    \begin{aligned}
    & \mathcal{V}(\tau_{m+1},\tilde{Z}(\tau_{m+1}),\tilde{Q}(\tau_{m+1})) - \mathcal{V}(t,\tilde{Z}(t),\tilde{Q}(t)) \\
    & \qquad \qquad = \int_{t}^{\tau_{m+1}} \nabla_z \mathcal{V}(s,\tilde{Z}(s),\tilde{Q}(s)) \cdot \dif W(s) - \int_{t}^{\tau_{m+1}} c \cdot \dif \tilde{L}(s).
    \end{aligned}
    \label{eq:equivalent_integral_finite_linear}
\end{equation}

Let $(\Omega, \mathcal{F}, (\mathcal{F}_t)_{t\in[0,T]}, \mathbb{P})$ denote the filtered probability space that supports the reference process and its driving diffusion process $W$. Processes $W$, $\tilde{Z}$ and $\tilde{Q}$ are all adapted to the filtration $(\mathcal{F}_t)_{t\in[0,T]}$. Then under suitable integrability assumptions, it follows that
\begin{equation}
\begin{aligned}
\E \Bigg[ \int_{t}^{\tau_{m+1}} \nabla_z \mathcal{V}(s,\tilde{Z}(s),\tilde{Q}(s)) \cdot \dif W(s) \Bigg| \mathcal{F}_t \Bigg] = 0, \quad \mathbb{P}\text{-a.s.}
\end{aligned}
\label{eq:zero_expectation}
\end{equation}
Combining \eqref{eq:equivalent_integral_finite_linear} and \eqref{eq:zero_expectation}, one arrives at the following for $m = 0, 1, \dots, M-1$ and $\mathbb{P}$-a.s.:
\begin{equation}
    \begin{aligned}
    & \E \Bigg[\mathcal{V}(\tau_{m+1},\tilde{Z}(\tau_{m+1}),\tilde{Q}(\tau_{m+1})) + \int_{t}^{\tau_{m+1}} c \cdot \dif \tilde{L}(s) \Bigg| \mathcal{F}_t \Bigg] = \mathcal{V}(t,\tilde{Z}(t),\tilde{Q}(t)) .
    \end{aligned}
\end{equation}
From the tower property of conditional expectations, we can write for $t \in (\tau_m, \tau_{m+1}]$ and $m = 0, 1, \dots, M-1$,
\begin{equation}
    \begin{aligned}
    & \E \Bigg[\mathcal{V}(\tau_{m+1},\tilde{Z}(\tau_{m+1}),\tilde{Q}(\tau_{m+1})) + \int_{t}^{\tau_{m+1}} c \cdot \dif \tilde{L}(s) \Bigg| \mathcal{F}_{\tau_m} \Bigg] =  \E \Bigg[\mathcal{V}(t,\tilde{Z}(t),\tilde{Q}(t)) \Bigg| \mathcal{F}_{\tau_m} \Bigg].
    \end{aligned}
    \label{eq:conditional_expectation_tau_m}
\end{equation}
In addition, letting $t\searrow \tau_m$ in \eqref{eq:V_linear_approx_eq}, we have for $m = 0, 1, \dots, M-1$,
\begin{equation}
    \begin{aligned}
    \lim_{t\searrow \tau_m} \mathcal{V}(t, z, q) = \mathcal{V}(\tau_m, z, q) - F(\tau_m, z, q, \mathcal{D} \mathcal{V}(\tau_m, z, q))  \Delta \tau_m.
    \end{aligned}
    \label{eq:V_linear_approx_eq_limit}
\end{equation}
Then, under suitable integrability assumptions, passing to the limit in \eqref{eq:conditional_expectation_tau_m} as $t \searrow \tau_m$ and substituting \eqref{eq:V_linear_approx_eq_limit} yields the following:
\begin{equation}
    \begin{aligned}
    & \E \Bigg[\mathcal{V}(\tau_{m+1},\tilde{Z}(\tau_{m+1}),\tilde{Q}(\tau_{m+1})) + \int_{\tau_m}^{\tau_{m+1}} c \cdot \dif \tilde{L}(s) \Bigg| \mathcal{F}_{\tau_m} \Bigg] \\
    & =  \lim_{t\searrow \tau_m} \E \Bigg[\mathcal{V}(t,\tilde{Z}(t),\tilde{Q}(t))\Bigg| \mathcal{F}_{\tau_m} \Bigg] = \E \Bigg[ \lim_{t\searrow \tau_m} \mathcal{V}(t,\tilde{Z}(t),\tilde{Q}(t))\Bigg| \mathcal{F}_{\tau_m} \Bigg]\\
    & = \mathcal{V}(\tau_m, \tilde{Z}(\tau_m),\tilde{Q}(\tau_m)) - F(\tau_m, \tilde{Z}(\tau_m), \tilde{Q}(\tau_m), \mathcal{D} \mathcal{V}(\tau_m, \tilde{Z}(\tau_m),\tilde{Q}(\tau_m)))  \Delta \tau_m.
    \end{aligned}
    \label{eq:conditional_expectation_tau_m_limit}
\end{equation}
Then taking the conditional expectations of both sides of \eqref{eq:conditional_expectation_tau_m_limit}  with respect to the $\sigma$-algebra generated by $\tilde{Z}(\tau_m)$ and $\tilde{Q}(\tau_m)$ and using the tower property yield
\begin{equation}
\begin{aligned}
& \E \Bigg[\mathcal{V}(\tau_{m+1},\tilde{Z}(\tau_{m+1}),\tilde{Q}(\tau_{m+1})) + c \cdot [\tilde{L}(\tau_{m+1}) - \tilde{L}(\tau_m)] \Bigg| \tilde{Z}(\tau_m), \tilde{Q}(\tau_m) \Bigg] \\
& = \mathcal{V}(\tau_m, \tilde{Z}(\tau_m),\tilde{Q}(\tau_m)) - F(\tau_m, \tilde{Z}(\tau_m), \tilde{Q}(\tau_m), \mathcal{D} \mathcal{V}(\tau_m, \tilde{Z}(\tau_m),\tilde{Q}(\tau_m)))  \Delta \tau_m.
\end{aligned}
\label{eq:conditional_expectation_tau_m_limit_final}
\end{equation}
Combining \eqref{eq:conditional_expectation_tau_m_limit_final} with \eqref{eq:V_approx} yields Equation \eqref{eq:deep_splitting_relationship_approximation}.

\section{Implementation Details}
\label{sec:implementation_details}

In this appendix, we record the numerical implementation used to produce the computational results in Section~\ref{sec:computational_results}. We describe the computational grid and network architecture, reference-path state sampling, two-stage backward training, demand generation, and policy evaluation. All implementation details correspond to the settings in Section~\ref{sec:data_description} unless otherwise stated.

\subsection{Computational Configuration}

The empirical baseline uses the weekly Fangonda instance described in Section~\ref{sec:data_description}, with $T=N=22$, $K=7$, and $h=1$. The deep splitting algorithm is run on a computational grid $0=\tau_0<\tau_1<\cdots<\tau_M=T$. The final policy is trained on a grid with $M=1100$, corresponding to 50 computational intervals per week. 

The baseline demand inputs are the calibrated weekly vectors $\mu(n)$ and covariance matrices $\Gamma(n)$ from Equations~\eqref{eq:weekly_demand_cha}--\eqref{eq:weekly_demand_hf}. On the finer computational grid, these weekly parameters are held constant over the corresponding week. The baseline cost vector is specified in Equation~\eqref{eq:shortage_cost_calibration}, where all cost quantities are measured in units of one untreated under-five malaria patient.

For each fine time interval $[\tau_m,\tau_{m+1}]$, the implementation uses a separate pair of feed-forward neural networks. The value network $H_m$ maps the extended state $(z,q)\in\mathbb R_+^{K+1}\times\mathbb R_+$ to a scalar approximation of $V(\tau_m,z,q)$. The gradient-difference network $G_m$ uses the same input and outputs a $(K+1)$-dimensional approximation of $\mathcal{D}V(\tau_m,z,q) = (\partial_{z_0} V - \partial_{q} V, \dots,\partial_{z_K}V - \partial_{q} V )$. For the Fangonda instance with $K=7$, both networks take a $9$-dimensional input; $H_m$ has a scalar output and $G_m$ has an $8$-dimensional output. Both $H_m$ and $G_m$ are fully connected feed-forward networks with GELU activation, 4 hidden layers, and 512 hidden neurons in each layer. The output layer of $H_m$ has no bias, while the output layer of $G_m$ includes a bias. All state coordinates are divided by the initial inventory $Q_0$ ($Q_0=142$ in the baseline) before entering the networks, and the output of $H_m$ is rescaled by a factor of $1/4$. Panel A of Table~\ref{tab:two_stage_training} records the network architecture parameters.

\subsection{Reference-Path State Sampling}

The networks are trained using the reference-path state-pool sampling method described in Section~\ref{sec:computational_method}. The reported baseline pool contains 32,768 paths: 16,384 from the fairness-only policy, 8,192 from the static policy, and 8,192 from the fluid-optimal policy.

Reference paths are initialized from a mixture of initial states. For $70\%$ of the paths, the initial condition is the physical initial state,
\[ \tilde{Z}(0)=0, \quad \tilde{Q}(0)=Q(0). \]
For the remaining $30\%$, the initial state is perturbed according to
\[ \tilde{Q}(0) = Q(0) + B_q, \quad B_q \sim \mathrm{Uniform}(-0.05Q(0),0.05Q(0)), \]
and
\[ \tilde{Z}_k(0)\sim \mathrm{Uniform}(0,0.0125Q(0)), \quad k=0,\ldots,K, \]
independently across locations. At each training time $\tau_m$, the states reached by these reference paths form a state pool $\mathcal{S}_m$. A training sample draws one complete state vector from $\mathcal{S}_m$; the components of $Z$ and $Q$ are not sampled independently. Starting from the sampled state, the implementation computes the reference allocation rate, projects it onto the one-step feasible set if needed, and advances the reflected Brownian reference process by one Euler--Maruyama step with a fresh Brownian increment.

\begin{figure}[!ht]
    \centering
    \includegraphics[width=\textwidth]{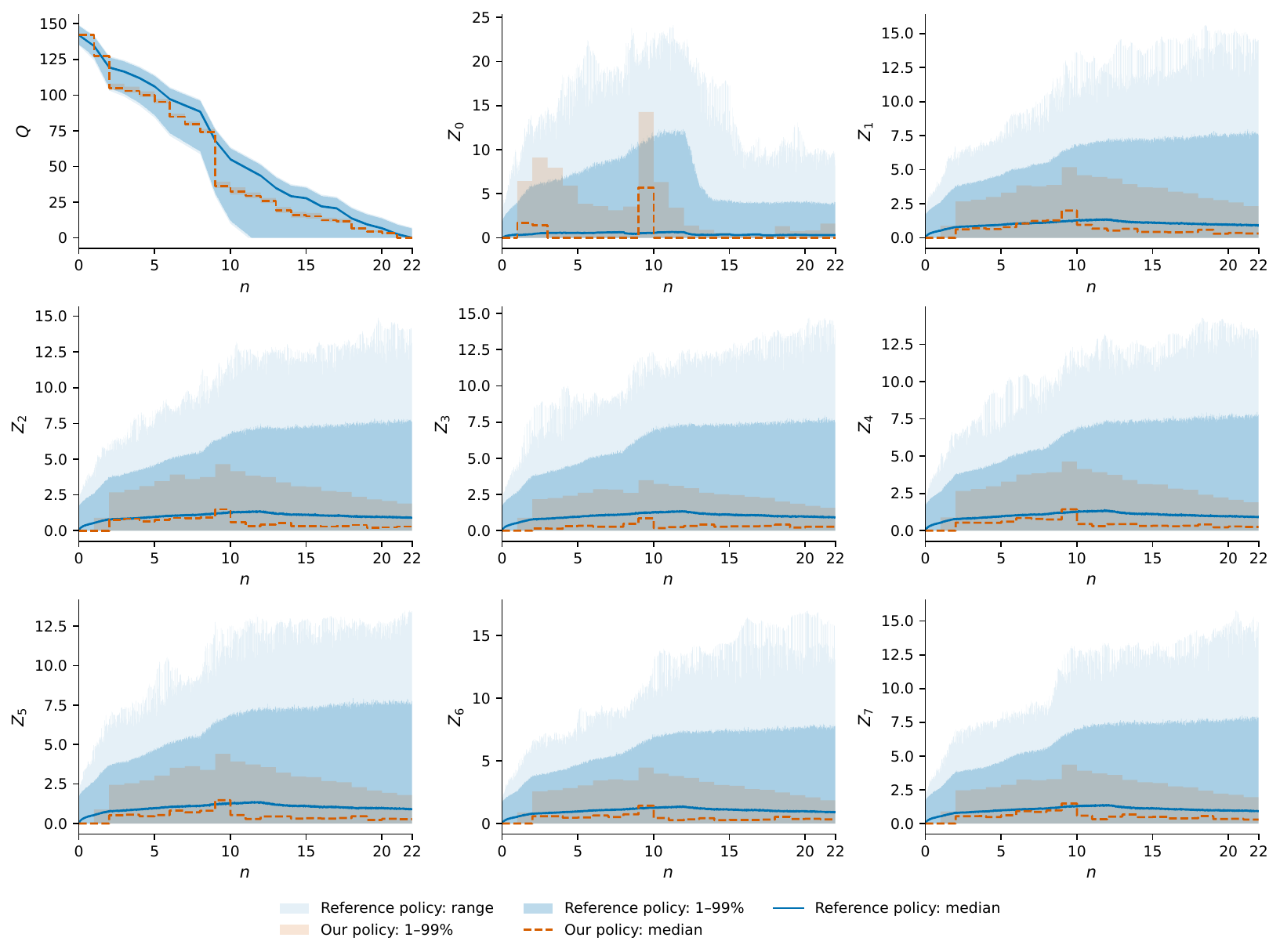}
    \caption{Componentwise state coverage in the baseline case.}
    \label{fig:baseline_state_coverage}
\end{figure}

We design the sampling method with a mixture of different reference policies and initial states so that the samples can better cover the states that our policy may visit, as Figure~\ref{fig:baseline_state_coverage} shows in the baseline case.

\subsection{Two-stage Backward Training}

Training proceeds backward over the fine computational grid $0=\tau_0<\tau_1<\cdots<\tau_M=T$. For the terminal interval, the regression target is the terminal condition. For each preceding interval, the already trained network at $\tau_{m+1}$ supplies the continuation target while the pair $(H_m,G_m)$ is optimized. 

The empirical loss is the sample-average version of the deep-splitting regression loss derived in Equations~\eqref{eq:empirical_loss_terminal}--\eqref{eq:empirical_loss_recursive}. In the implementation, we stabilize the joint training of the value network $H_m$ and the gradient-difference network $G_m$ by penalizing deviations between $G_m$ and the automatic-differentiation gradient difference of $H_m$. Specifically, the implemented training loss is
\[ \hat{\mathcal{J}}_m^{\mathrm{train}}(\omega_m,\nu_m) = \hat{\mathcal{J}}_m(\omega_m,\nu_m) + \zeta \hat C_m(\omega_m,\nu_m), \] where \[ \hat C_m(\omega_m,\nu_m) = \frac{1}{B(K+1)} \sum_{b=1}^B \left\| G_{\nu_m}^{(m)}(\tilde{Z}_m^{(b)},\tilde{Q}_m^{(b)}) - \mathcal D H_{\omega_m}^{(m)}(\tilde{Z}_m^{(b)},\tilde{Q}_m^{(b)}) \right\|_2^2. \]
In the reported runs, we use $\zeta=30$. We use a training batch size of 3,072 and a validation size of 16,384.

For each scenario and fairness weight, we use a two-stage backward training routine. We first train an auxiliary model using the computational method in Section~\ref{sec:computational_method} on a coarse grid with $M_c=110$ time intervals. For the terminal coarse interval, $j=M_c-1$, the value network and the gradient-difference network are randomly initialized using Xavier-uniform and Kaiming-uniform initialization, respectively. Training then proceeds backward, with each preceding network pair initialized from the trained pair at interval $j+1$.

\begin{figure}[!ht]
    \centering
    \includegraphics[width=0.75\linewidth]{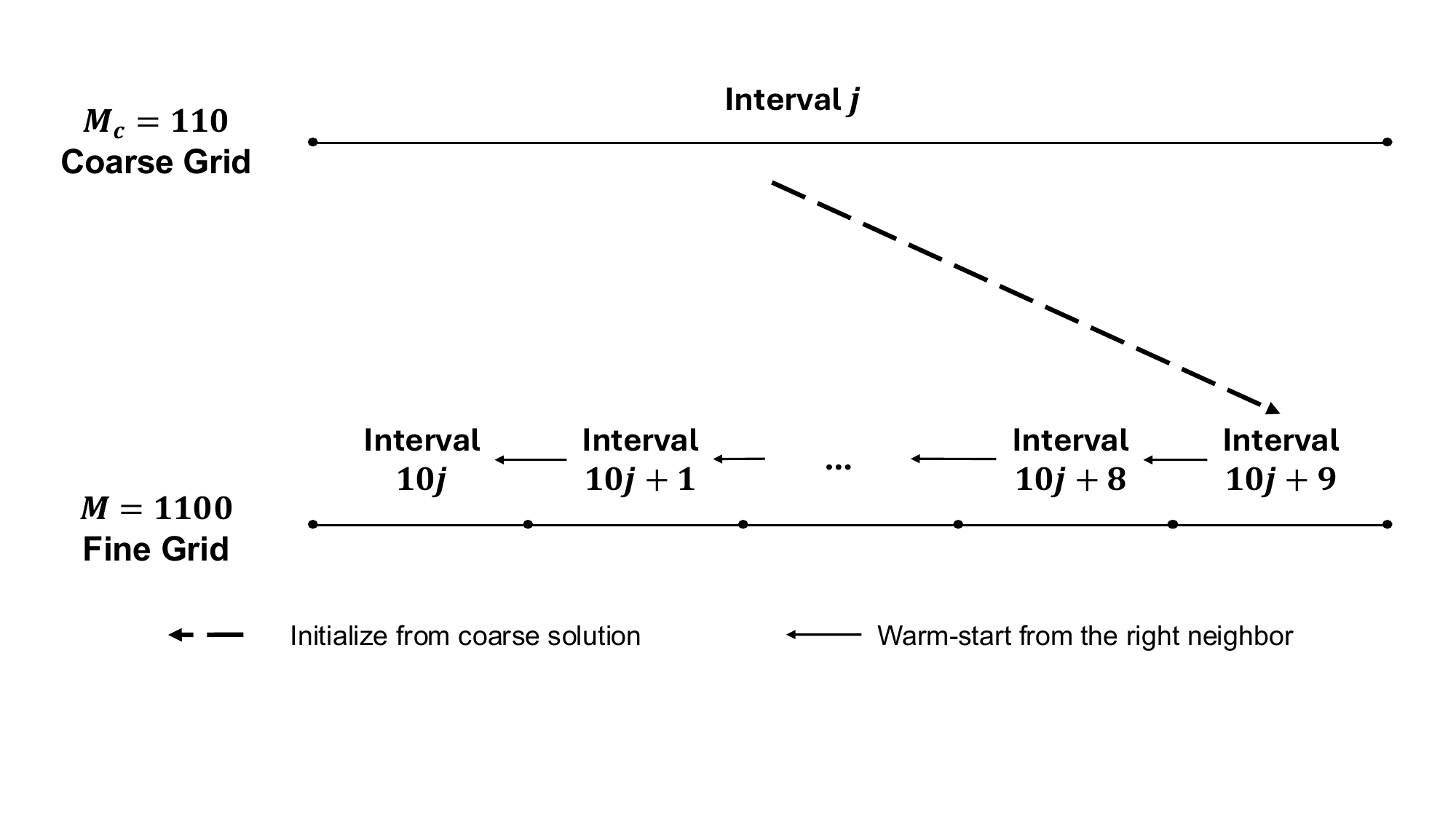}
    \caption{An illustration of blockwise initialization on the coarse interval $j$ and the corresponding ten fine intervals, $10j,10j+1,\dots,10j+9$.}
    \label{fig:block_init}
\end{figure}

We next train the model on a fine grid with $M=1{,}100$ time intervals using blockwise initialization from the coarse model. Specifically, coarse interval $j$ corresponds to fine intervals $10j,\ldots,10j+9$, for $j=0,\ldots,M_c-1$. The trained coarse network pair initializes the fine-grid pair at interval $10j+9$. Training then proceeds backward through that block, warm-starting each remaining pair from its trained fine-grid neighbor. Figure~\ref{fig:block_init} provides an illustration of this blockwise initialization.

For the reported experiments, the terminal interval of each stage uses the larger optimization budget without early stopping, whereas all nonterminal intervals use the smaller budget with early stopping. In particular, the terminal fine-grid interval is initialized from a coarse checkpoint but uses the terminal training schedule. The checkpoint with the lowest validation loss is restored for each interval. The training settings are summarized in Panels B--D of Table~\ref{tab:two_stage_training}.

\begin{table}[!ht]
    \centering
    \caption{Network architecture and training parameters
    for the two-stage backward training routine.}
    \label{tab:two_stage_training}
    \footnotesize
    \setlength{\tabcolsep}{5pt}
    \renewcommand{\arraystretch}{1.0}

    \begin{tabular}{lcc}
        \toprule
        Parameter & Coarse stage ($M_c=110$) & Fine stage ($M=1{,}100$) \\
        \midrule
        \multicolumn{3}{l}{\it Panel A: Network architecture.} \\
        Input dimension & \multicolumn{2}{c}{$9$} \\
        $H_m$/$G_m$ output dimension & \multicolumn{2}{c}{$1$, $8$} \\
        Number of hidden layers & \multicolumn{2}{c}{$4$} \\
        Number of neurons per layer & \multicolumn{2}{c}{$512$} \\
        Hidden-layer activation & \multicolumn{2}{c}{GELU} \\
        $H_m$/$G_m$ output bias & \multicolumn{2}{c}{No, Yes} \\ 
        Input normalization & \multicolumn{2}{c}{$(z,q)/Q_0$} \\
        $H_m$/$G_m$ output multiplier & \multicolumn{2}{c}{$0.25$, $1$} \\
        \midrule
        \multicolumn{3}{l}{\it Panel B: General training parameters.} \\
        Consistency coefficient $\zeta$ & \multicolumn{2}{c}{$30$} \\
        $H_m$/$G_m$ random initialization & \multicolumn{2}{c}{Xavier-uniform, Kaiming-uniform} \\
        Optimizer & \multicolumn{2}{c}{Adam} \\
        Training batch/validation sample size & \multicolumn{2}{c}{$3{,}072$, $16{,}384$} \\
        Joint gradient-clipping threshold & \multicolumn{2}{c}{$1$} \\
        \midrule
        \multicolumn{3}{l}{\it Panel C: Terminal interval training parameters.} \\
        Maximum iterations & $50{,}000$ & $30{,}000$ \\ 
        Learning rate & 
        \makecell{\scriptsize $10^{-3}$ $[0,10{,}000)$\\ \scriptsize $5\times 10^{-4}$ $[10{,}000,30{,}000)$ \\ \scriptsize $2.5\times 10^{-4}$ $[30{,}000,50{,}000)$} & \makecell{\scriptsize $10^{-3}$ $[0,6{,}000)$\\ \scriptsize $5\times 10^{-4}$ $[6{,}000,18{,}000)$ \\ \scriptsize $2.5\times 10^{-4}$ $[18{,}000,30{,}000)$} \\
        Early stopping & Disabled & Disabled \\
        \midrule
        \multicolumn{3}{l}{\it Panel D: Interior interval training parameters.}\\
        Maximum iterations & $8{,}000$ & $2{,}500$ \\
        Learning rate & 
        \makecell{\scriptsize $2\times10^{-4}$ $[0,1{,}000)$\\ \scriptsize $10^{-4}$ $[1{,}000,4{,}000)$ \\ \scriptsize $5\times 10^{-5}$ $[4{,}000,8{,}000)$} & \makecell{\scriptsize $2\times10^{-4}$ $[0,400)$\\ \scriptsize $10^{-4}$ $[400,1{,}200)$ \\ \scriptsize $5\times 10^{-5}$ $[1{,}200,2{,}500)$} \\
        Early stopping & Enabled & Enabled \\
        Minimum iterations & $1{,}200$ & $1{,}200$ \\
        Patience (iterations) & $1{,}500$ & $1{,}500$ \\
        Minimum validation improvement & $5\times10^{-5}$ & $5\times10^{-5}$ \\
        \bottomrule
    \end{tabular}
\end{table}

\subsection{Demand Generation}
\label{app:demand_generation}

Recall that the weekly demand is generated from the calibrated binomial model:
\[ \xi_0(n)\sim \mathrm{Binomial}(N_h,\theta(n)), \quad \xi_k(n)\sim \mathrm{Binomial}(N_c,\theta(n)), \ k=1,\ldots,K. \]

\paragraph{CHA demand variability.}
For variance multiplier $c\in\{2,5\}$, set $m_c=(N_c-c)/(c-1)$. Independently for each CHA $k$ and week $n$, we first draw $P_k(n)$ and then draw $\xi_k(n)$ conditional on \(P_k(n)\) as follows:
\[ P_k(n)\sim\operatorname{Beta}(m_c\theta(n),m_c[1-\theta(n)]), \qquad \xi_k(n)\mid P_k(n)\sim\operatorname{Binomial}(N_c,P_k(n)). \]
The marginal distribution of CHA demand has mean $N_c\,\theta(n)$ and variance $cN_c\,\theta(n)[1-\theta(n)]$. The case $c=1$ is the baseline binomial model. HF demand remains binomial with its baseline mean and variance. Shipment bounds and the demand means are unchanged.

\paragraph{Spatial demand correlation.}
In this experiment, we introduce correlation between the weekly demands at different locations while preserving their baseline binomial distributions. The demands considered are local HF demand, $\xi_0(n)$, and CHA demands, $\xi_k(n)$, $k=1,\ldots,K$. We consider a target pairwise Pearson correlation of $\rho\in\{0.5,0.7\}$ within each week. Demand vectors remain independent across weeks.

We generate correlated demand using a Gaussian copula. Let $F_{k,n}$ denote the baseline binomial cumulative distribution function for demand at location $k$ in week $n$. Thus, $F_{0,n}$ corresponds to $\operatorname{Binomial}(N_h,\theta(n))$, whereas $F_{k,n}$ corresponds to $\operatorname{Binomial}(N_c,\theta(n))$ for $k=1,\ldots,K$. Let $\Phi$ denote the standard normal cumulative distribution function, and let $F_{k,n}^{-1}$ denote the binomial quantile function. If $Y$ is standard normal, then
\[ F_{k,n}^{-1}(\Phi(Y)) \]
has the required binomial distribution. We therefore construct correlated binomial demands by applying this transformation to correlated standard normal variables.

We first calibrate the correlations of these normal variables. The transformation generally changes Pearson correlations, so the normal correlation must be chosen to produce the desired demand correlation $\rho$. For each week $n$ and pair of locations $0\leq k<\ell\leq K$, we choose $c_{k\ell,n}$ such that
\[ \operatorname{Corr}\left( F_{k,n}^{-1}(\Phi(U)), F_{\ell,n}^{-1}(\Phi(V)) \right)=\rho, \]
where $(U,V)$ is a bivariate standard normal vector with correlation $c_{k\ell,n}$. We solve this equation by bisection, evaluating the demand correlation using binomial cumulative probabilities and the bivariate normal cumulative distribution function.

We next assemble the calibrated correlations into a symmetric matrix with diagonal entries equal to one. Because the pairwise calibration does not guarantee that this matrix is positive definite, we apply an iterative adjustment: we replace eigenvalues below $10^{-10}$ by $10^{-10}$, rescale the matrix to restore unit diagonal entries, and restrict its off-diagonal entries to $[-0.9999,0.9999]$. Let $C_n$ denote the resulting correlation matrix. This adjustment can change the implied demand correlations, so we check the correlations after the adjustment; they may differ from the target $\rho$.

Having constructed $C_n$, we generate demand for week $n$ by drawing
\[ Y(n)\sim\mathcal N(0,C_n) \]
and setting
\[ \xi_k(n)=F_{k,n}^{-1}(\Phi(Y_k(n))), \qquad k=0,\ldots,K. \]
This construction preserves each location's baseline binomial distribution and hence its mean and variance. We draw the vectors $Y(n)$ independently across weeks. The case $\rho=0$ corresponds to the independent baseline and uses $C_n=I$.
Finally, for the Brownian approximation used to train the policy, we retain the baseline marginal variances and set
\[ \Gamma_{k\ell}(n) =(C_n)_{k\ell} \sqrt{\Gamma_{kk}(n)\Gamma_{\ell\ell}(n)}, \qquad k\ne\ell. \]
Thus, the Brownian approximation uses the correlations of the auxiliary normal variables. Its covariance matrix generally differs from that of the binomial demands generated above.

\subsection{Policy Evaluation}

For each experimental setting, our policy and benchmark policies are evaluated in the original discrete-time model of Section~\ref{sec:model} using $10{,}000$ Monte Carlo demand paths. Each path starts from $Z(0)=0$ and the initial warehouse inventory $Q(0)$ specified for that setting. In each period $n$, the number of referrals from CHA $k$ to the HF is set to $p_k[L_k(n)-L_k(n-1)]$. All policies are evaluated on the same demand paths using common random numbers, with the trained neural networks held fixed throughout evaluation. 

For each policy and sample path, we record the objective decomposition and operational metrics defined in Appendix~\ref{app:evaluation_metrics}. We compute Monte Carlo means and standard errors across paths using the formulas given there. Unless otherwise stated, the tables report Monte Carlo means. Cost quantities are expressed in units of one untreated under-five malaria patient.

\end{document}

%% file: math_command.tex
\newcommand{\E}{\mathbb{E}}					%
\newcommand{\R}{\mathbb{R}}					%
\DeclareMathOperator*{\argmin}{argmin}

\renewcommand{\geq}{\geqslant}
\renewcommand{\leq}{\leqslant}

\newcommand{\dif}{\mathop{}\!\mathrm{d}} %